%% file: agt-3-10.tex
\documentclass{gtart}

\input agtout

\lognumber{10}
\volumenumber{3}
\volumeyear{2003}
\papernumber{10}
\published{5 March 2003}
\pagenumbers{235}{285}
\received{12 November 2001}
\revised{5 February 2003}
\accepted{14 February 2003}

\usepackage{amsmath,amssymb,graphicx}

\newcommand{\face}[1]{\xrightarrow{#1}}
\newcommand{\inv}{^{-1}}

\newcommand{\bN}{\mathbf{N}}
\newcommand{\bR}{\mathbf{R}}
\newcommand{\bZ}{\mathbf{Z}}
\newcommand{\za}{\alpha}
\newcommand{\zb}{\beta}
\newcommand{\zd}{\delta}
\newcommand{\ze}{\epsilon}
\newcommand{\zg}{\gamma}
\newcommand{\zj}{\psi}

\newcommand{\zr}{\rho}
\newcommand{\zs}{\sigma}
\newcommand{\zt}{\tau}
\newcommand{\zv}{\varphi}
\newcommand{\zS}{\Sigma}

\newtheorem{thm}{Theorem}
\theoremstyle{definition}
\newtheorem{ex}[thm]{Example}

\newcommand{\sections}{\renewcommand{\thethm}{\thesection.\arabic{thm}}
           \setcounter{thm}{0}}
\newcommand{\subsections}{\setcounter{thm}{0}\renewcommand{\thethm}
           {\thesubsection.\arabic{thm}}}
\newcommand{\nosubsections}{\renewcommand{\thethm}{\thesection.\arabic{thm}}
           \setcounter{thm}{0}}
\newcommand{\nosubsubsections}{\setcounter{thm}{0}
           \renewcommand{\thethm}{\thesubsection.\arabic{thm}}}
\newcommand{\linnum}{\stepcounter{thm}\tag{\thethm}}

\begin{document}
\title[Heegaard diagrams and surgery descriptions]{Heegaard
diagrams and surgery descriptions\\for twisted face-pairing 3-manifolds}
\author{J.W. Cannon\\W.J. Floyd\\W.R. Parry}
\shortauthors{Cannon, Floyd and Parry}

\address{Department of Mathematics, Brigham Young University\\Provo, UT
84602, USA\\\smallskip\\Department of Mathematics, Virginia 
Tech\\Blacksburg, VA 24061, USA\\\smallskip\\Department of 
Mathematics, Eastern Michigan University\\Ypsilanti, MI 48197, USA} 
\asciiaddress{Department of Mathematics, Brigham Young University\\Provo, UT
84602, USA\\\\Department of Mathematics, Virginia 
Tech, Blacksburg, VA 24061, USA\\\\Department of 
Mathematics, Eastern Michigan University\\Ypsilanti, MI 48197, USA} 
\email{cannon@math.byu.edu, floyd@math.vt.edu,
walter.parry@emich.edu}
\url{http://www.math.vt.edu/people/floyd}

\begin{abstract}
The twisted face-pairing construction of our earlier papers gives an
efficient way of generating, mechanically and with little effort, myriads
of relatively simple face-pairing descriptions of interesting closed
3-manifolds. The corresponding description in terms of surgery, or
Dehn-filling, reveals the twist construction as a carefully organized
surgery on a link.

In this paper, we work out the relationship between the twisted
face-pairing description of closed 3-manifolds and the more common
descriptions by surgery and Heegaard diagrams. We show that all Heegaard
diagrams have a natural decomposition into subdiagrams called Heegaard
cylinders, each of which has a natural shape given by the ratio of two
positive integers. We characterize the Heegaard diagrams arising naturally
from a twisted face-pairing description as those whose Heegaard cylinders
all have integral shape. This characterization allows us to use the Kirby
calculus and standard tools of Heegaard theory to attack the problem of
finding which closed, orientable 3-manifolds have a twisted face-pairing
description.
\end{abstract}
\primaryclass{57N10} \keywords{3-manifold constructions, Dehn surgery,
Heegaard diagrams}
\maketitle

\sections

\section{Introduction}\label{sec:intro}\nosubsections

The twisted face-pairing construction of our earlier papers
\cite{intro}, \cite{twisted}, \cite{ample} gives an efficient way of
generating, mechanically and with little effort, myriads of
relatively simple face-pairing descriptions of interesting closed
3-manifolds. Starting with a faceted 3-ball $P$ and an arbitrary
orientation-reversing face-pairing $\epsilon$ on $P$, one constructs
a faceted 3-ball $Q$ and an orientation-reversing face-pairing
$\delta$ on $Q$ such that the quotient $Q/\delta$ is a manifold. Here
$Q$ is obtained from $P$ by subdividing the edges according to a
function which assigns a positive integer (called a multiplier) to
each edge cycle, and $\delta$ is obtained from $\epsilon$ by
precomposing each face-pairing map with a twist. Which direction to
twist depends on the choice of an orientation of $P$. Hence for a
given faceted 3-ball $P$, orientation-reversing face-pairing
$\epsilon$, and multiplier function, one obtains two twisted
face-pairing manifolds $M=Q/\delta$ and $M^*=Q/\delta^*$ (one for
each orientation of $P$).

In \cite{intro} and \cite{twisted} we introduced twisted face-pairing
3-manifolds  and developed their first properties. A surprising result in
\cite{twisted} is the duality theorem that says that, if $P$ is a regular
faceted 3-ball, then $M$ and $M^*$ are homeomorphic
in a way that makes their cell structures dual to each other. This duality is
instrumental in \cite{ample}, where we investigated a special subset of
these manifolds, the ample twisted face-pairing manifolds.  We showed that
the fundamental group of every ample twisted face-pairing manifold is
Gromov hyperbolic with space at infinity a 2-sphere.

In this paper we connect the twisted face-pairing construction with two
standard 3-manifold constructions. Starting with a faceted 3-ball $P$ with
$2g$ faces and an orientation-reversing face-pairing $\ze$ on $P$, we
construct a closed surface $S$ of genus $g$ and two families $\zg$ and
$\zb$ of pairwise disjoint simple closed curves on $S$. The elements of
$\zg$ correspond to the face pairs and the elements of $\zb$ correspond to
the edge cycles of $\ze$. Given a choice of multipliers for the edge
cycles, we then give a Heegaard diagram for the resulting twisted
face-pairing 3-manifold. The surface $S$ is the Heegaard surface, and the
family $\zg$ is one of the two families of meridian curves. The other
family is obtained from $\zg$ by a product of powers of Dehn twists along
elements of $\zb$; the powers of the Dehn twists are the multipliers. From
the Heegaard diagram, one can easily construct a framed link in the
3-sphere such that Dehn surgery on this framed link gives the twisted
face-pairing manifold. The components of the framed link fall naturally
into two families; each curve in one family corresponds to a face pair and
has framing 0, and each curve in the other family corresponds to an edge
cycle and has framing the sum of the reciprocal of its multiplier and the
blackboard framing of a certain projection of the curve. These results are
very useful for understanding both specific face-pairing manifolds and
entire classes of examples. While we defer most illustrations of these
results to a later paper \cite{survey}, we give several examples here to
illustrate how to use these results to give familiar names to some twisted
face-pairing 3-manifolds.

\begin{figure}[ht!]
\centerline{\includegraphics{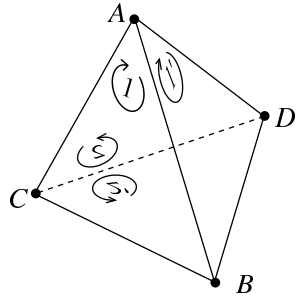}} \caption{The complex $P$}
\label{fig:tet1p}
\end{figure}

One of our most interesting results in this paper is that all Heegaard
diagrams have a natural decomposition into subdiagrams called Heegaard
cylinders, each of which has a natural shape given by the ratio of two
positive integers. We characterize the Heegaard diagrams arising naturally
from a twisted face-pairing description as those whose Heegaard cylinders
all have integral shape.

We give a preliminary example to illustrate the twisted face-pairing
construction. Let $P$ be a tetrahedron with vertices $A$, $B$, $C$, and
$D$, as shown in Figure~\ref{fig:tet1p}. Consider the face-pairing
$\ze=\{\ze_1,\ze_2\}$ on $P$ with map $\epsilon_1$ which takes triangle
$ABC$ to triangle $ABD$ fixing the edge $AB$ and map $\epsilon_2$ which
takes triangle $ACD$ to $BCD$ fixing the edge $CD$. This example was
considered briefly in \cite{intro} and in more detail in \cite[Example
3.2]{twisted}. The edge cycles are the equivalence classes of the edges of
$P$ under the face-pairing maps. The three edge cycles are $\{AB\}$,
$\{BC,BD,AD,AC\}$, and $\{CD\}$; the associated diagrams of face-pairing
maps are shown below.
\begin{equation*}
AB\face{\ze_1}AB
\end{equation*}
\begin{equation*}
BC\face{\ze_1}BD\face{\ze_2\inv}AD\face{\ze_1\inv}AC\face{\ze_2}BC
\end{equation*}
\begin{equation*}
CD\face{\ze_2}CD
\end{equation*}
To construct a twisted face-pairing manifold from $P$, for each edge cycle
$[e]$ we choose a positive integer $\text{mul}([e])$ called the multiplier
of $[e]$. Let $Q$ be the subdivision of $P$ obtained by subdividing each
edge $e$ of $P$ into $\#([e])\cdot \text{mul}([e])$ subedges. The
face-pairing maps $\ze_1$ and $\ze_2$ naturally give face-pairing maps on
the faces of $Q$. Choose an orientation of $\partial Q$, and
define the twisted face-pairing $\zd$ on $Q$ by precomposing each $\ze_i$
with an orientation-preserving homeomorphism of its domain which takes
each vertex to the vertex that follows it in the induced orientation on
the boundary. By the fundamental theorem of twisted face-pairings (see
\cite{intro} or \cite{twisted}), the quotient $Q/\zd$ is a closed
3-manifold.

\begin{figure}[ht!]
\centerline{\includegraphics{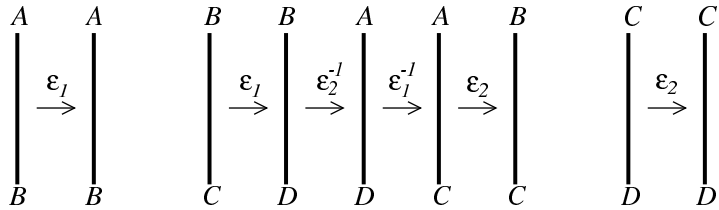}} \caption{The edge diagrams}
\label{fig:tet1pdia}
\end{figure}

\begin{figure}[ht!]
\centerline{\includegraphics{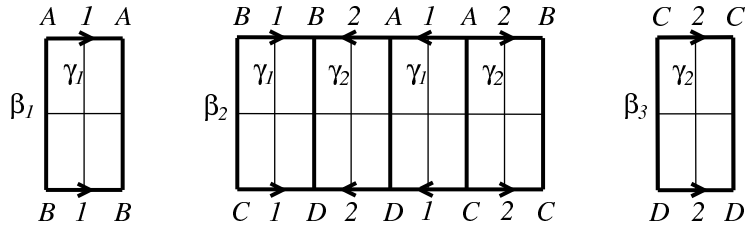}} \caption{The rectangles that
correspond to the edge diagrams} \label{fig:tet1pdi2}
\end{figure}

\begin{figure}[ht!]
\centerline{\includegraphics{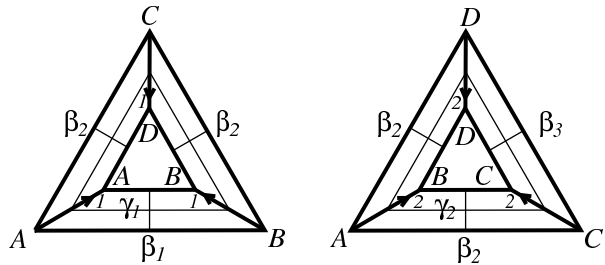}} \caption{Another view of the
surface $S$} \label{fig:tet1pdi4}
\end{figure}

\begin{figure}[ht!]
\centerline{\includegraphics{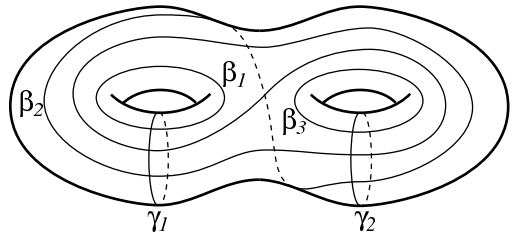}} \caption{The curve families
$\{\gamma_1,\gamma_2\}$ and $\{\beta_1,\beta_2,\beta_3\}$ on the surface
$S$} \label{fig:tet1phe}
\end{figure}

\begin{figure}[ht!]
\centerline{\includegraphics{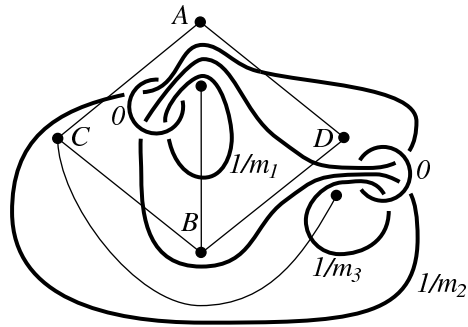}} \caption{The surgery
description} \label{fig:tet1pdi3}
\end{figure}

To construct a Heegaard diagram and framed link for the twisted
face-pairing manifold $Q/\zd$, we first schematically indicate the edge
diagrams as shown in Figure~\ref{fig:tet1pdia}. We then make rectangles
out of the edge diagrams in Figure~\ref{fig:tet1pdi2}, and add thin
horizontal and vertical line segments through the midpoints of each of the
subrectangles of the rectangles. We identify the boundary edges of the
rectangles in pairs preserving the vertex labels (and, for horizontal
edges, the order) to get a quotient surface $S$ of genus two. The image in
$S$ of the thin vertical arcs is a union of two disjoint simple closed
curves $\zg_1$ and $\zg_2$, which correspond to the two face pairs. The
image in $S$ of the thin horizontal arcs is a union of three pairwise
disjoint simple closed curves $\zb_1$, $\zb_2$, and $\zb_3$, which
correspond to the three edge cycles. Figure~\ref{fig:tet1pdi4} shows $S$
as the quotient of the union of two annuli, and Figure~\ref{fig:tet1phe}
shows the curve families $\{\zg_1,\zg_2\}$ and $\{\zb_1,\zb_2,\zb_3\}$ on
$S$. For $i\in\{1,2,3\}$, let $m_i$ be the multiplier of the edge cycle
corresponding to $\zb_i$ and let $\tau_i$ be one of the two Dehn twists
along $\zb_i$. We choose $\zt_1$, $\zt_2$, and $\zt_3$ so that they are
oriented consistently. Let $\zt = \zt_1^{m_1} \circ \zt_2^{m_2} \circ
\zt_3^{m_3}$. It follows from Theorem~\ref{thm:dehntwist} that $S$ and
$\{\zg_1,\zg_2\}$ and $\{\zt(\zg_1),\zt(\zg_2)\}$ form a Heegaard diagram
for the twisted face-pairing manifold $Q/\zd$. From the Heegaard diagram,
one can use standard techniques to give a framed surgery description for
$Q/\zd$. An algorithmic description for this is given in
Theorem~\ref{thm:surgery}. In the present example, the surgery description
is shown in Figure~\ref{fig:tet1pdi3} together with a modification of the
1-skeleton of the tetrahedron $P$. There are two curves with framing $0$,
corresponding to the two pairs of faces. The other three curves correspond
to the edge cycles and have framings the reciprocals of the multipliers.

We now describe our Heegaard diagram construction in greater detail.  We
use the notation and terminology of \cite{twisted}.  Let $P$ be a faceted
3-ball, let $\ze$ be an orientation-reversing face-pairing on $P$, and let
mul be a multiplier function for $\ze$. (As in \cite{twisted}, we for now
assume that $P$ is a regular CW complex. We drop the regularity assumption
in Section~\ref{sec:generalizing}.) Let $Q$ be the twisted face-pairing
subdivision of $P$, let $\zd$ be the twisted face-pairing on $Q$, and let
$M$ be the associated twisted face-pairing manifold.  We next construct a
closed surface $S$ with the structure of a cell complex.  For this we fix
a cell complex $X$ cellularly homeomorphic to the 1-skeleton of $Q$.
Suppose given two paired faces $f$ and $f^{-1}$ of $Q$.  We choose one of
these faces, say $f$, and we construct $\partial f\times[0,1]$.  We view
the interval $[0,1]$ as a 1-cell, and we view $\partial f\times[0,1]$ as a
2-complex with the product cell structure. For every $x\in\partial f$ we
identify $(x,0)\in\partial f\times[0,1]$ with the point in $X$
corresponding to $x$ and we identify $(x,1)\in\partial f\times[0,1]$ with
the point in $X$ corresponding to $\zd_f(x)\in\partial f^{-1}$.  Doing
this for every pair of faces of $Q$ yields a cell complex $Y$ on a closed
surface.  We define $S$ to be the first dual cap subdivision of $Y$;
because every face of $Y$ is a quadrilateral, this simply means that to
obtain $S$ from $Y$ we subdivide every face of $Y$ into four
quadrilaterals in the straightforward way.  We say that an edge of $S$ is
vertical if it is either contained in $X$ or is disjoint from $X$.  We say
that an edge of $S$ is diagonal if it is not vertical.  The union of the
vertical edges of $S$ which are not edges of $Y$ is a family of simple
closed curves in $S$.  Likewise the union of the diagonal edges of $S$
which are not edges of $Y$ is a family of simple closed curves in $S$.
Theorem~\ref{thm:twistedheegaard} states that the surface $S$ and these
two families of curves form a Heegaard diagram for $M$.

\begin{figure}[ht!]
\centerline{\includegraphics{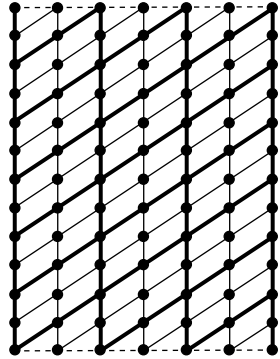}} \caption{The cylinder $C$
corresponding to the edge cycle $E$} \label{fig:cylinder}
\end{figure}

In this paragraph we indicate how to associate to a given edge cycle $E$
of $\ze$ a closed subspace of $S$.  To simplify this discussion we assume
that $E$ contains three edges and that $\text{mul}(E)=2$.  When
constructing $Q$ from $P$, every edge of $E$ is subdivided into $2\cdot
3=6$ subedges.  So corresponding to the three edges of $E$, the complex
$S$ contains three 1-complexes, each of them homeomorphic to an interval
and the union of 12 vertical edges of $S$.  These three 1-complexes and
part of $S$ are shown in Figure~\ref{fig:cylinder}; the three 1-complexes
are drawn as four thick vertical line segments with the left one to be
identified with the right one.  We refer to the closed subspace $C$ of $S$
shown in Figure~\ref{fig:cylinder} as an edge cycle cylinder or simply as
a cylinder.  In Figure~\ref{fig:cylinder}, vertical edges of $S$ are drawn
vertically and diagonal edges of $S$ are drawn diagonally.  Some arcs in
Figure~\ref{fig:cylinder} are dashed because they are not contained in the
1-skeleton of $S$.  The thick edges in Figure~\ref{fig:cylinder} are the
edges of $Y$ in $C$.  (It is interesting to note that these thick edges
essentially give the diagram in Figure 11 of \cite{twisted}.)  Note that
the edge cycle cylinder $C$ need not be a closed annulus, although $C$ is
the closure of an open annulus. (Identifications of boundary points are
possible.)  We choose these edge cycle cylinders so that their union is
$S$ and the cylinders of distinct $\ze$-edge cycles have disjoint
interiors.

We define the circumference of an edge cycle cylinder to be the number of
edges in its edge cycle.  We define the height of an edge cycle cylinder
to be the number of edges in its edge cycle times the multiplier of its
edge cycle.  The edge cycle cylinder $C$ in Figure~\ref{fig:cylinder}
contains three arcs $\zr_1$, $\zr_2$, $\zr_3$ whose endpoints lie on
dashed arcs such that each of $\zr_1$, $\zr_2$, $\zr_3$ is a union of thin
vertical edges.  Likewise $C$ contains three arcs $\zs_1$, $\zs_2$,
$\zs_3$ such that each of $\zs_1$, $\zs_2$, $\zs_3$ is a union of thin
diagonal edges and the endpoints of $\zs_i$ equal the endpoints of $\zr_i$
for every $i\in\{1,2,3\}$.  Because the height of $C$ equals 2 times the
circumference of $C$, it follows that $\zs_1$, $\zs_2$, $\zs_3$ can be
realized as the images of $\zr_1$, $\zr_2$, $\zr_3$ under the second power
of a Dehn twist along a waist of $C$.  This observation and the previous
paragraphs essentially give the following.  Let $\za_1,\dotsc,\za_n$ be
the simple closed curves in $S$ which are unions of vertical edges of $S$
that are not edges of $Y$.  Let $E_1,\dotsc,E_m$ be the edge cycles of
$\ze$. For every $i\in\{1,\dotsc,m\}$ construct a waist $\zb_i$ in the
edge cycle cylinder of $E_i$ so that $\zb_1,\dotsc,\zb_m$ are pairwise
disjoint simple closed curves in $S$.  For every $i\in\{1,\dotsc,m\}$ let
$\zt_i$ be one of the two Dehn twists on $S$ along $\zb_i$, chosen so that
the directions in which we twist are consistent.  Set
$\zt^{\text{mul}}=\zt_1^{\text{mul}(E_1)}\circ\cdots\circ
\zt_m^{\text{mul}(E_m)}$.  Then $S$ and $\za_1,\dotsc,\za_n$ and
$\zt^{\text{mul}}(\za_1),\dotsc,\zt^{\text{mul}}(\za_n)$ form a Heegaard
diagram for $M$.  The last statement is the content of
Theorem~\ref{thm:dehntwist}.

The result of the previous paragraph leads to a link $L$ in $S^3$
that has components $\zg_1,\dotsc,\zg_n$ and
$\zd_1,\dotsc,\zd_m$, where $\zg_1,\dotsc,\zg_n$ correspond to
$\za_1,\dotsc,\za_n$ and $\zd_1,\dotsc,\zd_m$ correspond to
$\zb_1,\dotsc,\zb_m$.  We define a framing of $L$ so that
$\zg_1,\dotsc,\zg_n$ have framing 0 and for every
$i\in\{1,\dotsc,m\}$ $\zd_i$ has framing $\text{mul}(E_i)^{-1}$ plus
the blackboard framing of $\zd_i$ relative to a certain projection.
Then the manifold obtained by Dehn surgery on $L$ is homeomorphic to
$M$.  The last statement is the content of Theorem~\ref{thm:surgery}.
At last we see that multipliers of edge cycles are essentially
inverses of framings of link components.  In
Section~\ref{sec:corridor} we make the construction of $L$
algorithmic and simple using what we call the corridor construction.

We know of no nice characterization of twisted face-pairing
3-manifolds. 
However, Theorem~\ref{thm:cylinder} gives such a characterization of
their Heegaard diagrams.  Theorem~\ref{thm:cylinder} and results leading
to it give the following statements.
Every irreducible Heegaard diagram for an
orientable closed 3-manifold $M$ gives rise to a faceted 3-ball $P$ with
orientation-reversing face-pairing $\ze$ (in essentially two ways -- one
for each family of meridian curves) such that $P/\ze $ is homeomorphic
to $M$.
Every irreducible Heegaard diagram can be decomposed into
cylinders, which we call Heegaard cylinders, essentially just as our
above Heegaard diagrams of twisted face-pairing manifolds are decomposed
into edge cycle cylinders.
In general heights of Heegaard cylinders are not multiples of
their circumferences.  A given irreducible Heegaard diagram is the
Heegaard diagram, as constructed above, of a twisted face-pairing
manifold if and only if the height of each of its Heegaard cylinders is
a multiple of its circumference.  Furthermore, if the height of every
Heegaard cylinder is a multiple of its circumference, then the
face-pairing $\ze$ constructed from the given Heegaard diagram is a
twisted face-pairing.

Thus far we have discussed the construction of Heegaard diagrams for
twisted face-pairing manifolds and the construction of face-pairings
from irreducible Heegaard diagrams.  In Theorem~\ref{thm:heegaard} we more
generally construct (irreducible) Heegaard diagrams for manifolds of the form
$P/\ze $, where $P$ is a faceted 3-ball with orientation-reversing
face-pairing $\ze$ and the cell complex $P/\ze $ is a manifold with one
vertex.  In Theorem~\ref{thm:cylinder} we construct for every
irreducible Heegaard diagram for a 3-manifold $M$ a faceted 3-ball $P$
with orientation-reversing face-pairing $\ze$ (in essentially two ways
-- one for each family of meridian curves) such that $P/\ze $ is a cell
complex with one vertex and $P/\ze $ is homeomorphic to $M$.  These two
constructions are essentially inverse to each other.

The above statements that every irreducible Heegaard diagram gives rise to
a faceted 3-ball require a more general definition of faceted 3-ball than
the one given in \cite{twisted}.  In \cite{twisted} faceted 3-balls are
regular, that is, for every open cell of a faceted 3-ball the prescribed
homeomorphism of an open Euclidean ball to that cell extends to a
homeomorphism of the closed Euclidean ball to the closed cell.  On the
other hand, the cellulation of the boundary of a 3-ball which arises from
a Heegaard diagram has paired faces but otherwise is arbitrary.  So we now
define a faceted 3-ball $P$ to be an oriented CW complex such that $P$ is
a closed 3-ball, the interior of $P$ is the unique open 3-cell of $P$, and
the cell structure of $\partial P$ does not consist of just one 0-cell and
one 2-cell.  This generalization presents troublesome minor technical
difficulties but no essential difficulties. In particular, all the results
of \cite{intro} and \cite{twisted} hold for these more general faceted
3-balls.  Section~\ref{sec:generalizing} deals with this generalization.
Except when the old definition is explicitly discussed, we henceforth in
this paper use the new definition of faceted 3-ball.  We know of no
reducible twisted face-pairing manifold which arises from a regular
faceted 3-ball; the old twisted face-pairing manifolds all seem to be
irreducible.  On the other hand the new twisted face-pairing manifolds are
often reducible.  See Examples~\ref{ex:target}
(which is considered again in \ref{ex:target'} and \ref{ex:target''})
and \ref{ex:s2xs1} (which is considered again in \ref{ex:s2xs1'}).

Our construction of Heegaard diagrams from face-pairings uses a
subdivision of cell complexes which we call dual cap subdivision.  We
define and discuss dual cap subdivision in Section~\ref{sec:dualcap}. The
term ``dual'' is motivated by the notion of dual cell complex, and the
term ``cap'' is motivated by its association with intersection.
Intuitively, the dual cap subdivision of a cell complex is gotten by
``intersecting'' the complex with its ``dual complex''.  Dual cap
subdivision is coarser than barycentric subdivision, and it is well suited
to the constructions at hand.  Heegaard decompositions of 3-manifolds are
usually constructed by triangulating the manifolds and working with their
second barycentric subdivisions.  Instead of using barycentric
subdivision, we use dual cap subdivision, and we obtain the following.
Earlier in the introduction we constructed a surface $S$ with a
cell structure.
We show that $S$ is cellularly homeomorphic to
a subcomplex of the second dual cap subdivision of the manifold $M$, where
this subcomplex corresponds to the usual Heegaard surface gotten by using
a triangulation and barycentric subdivision.

In Section~\ref{sec:examples} we use the corridor construction of
Section~\ref{sec:corridor} to construct links in $S^3$ for three
different model face-pairings.  Simplifying these links using isotopies
and Kirby calculus, we are able to identify the corresponding twisted
face-pairing manifolds.  In Example~\ref{ex:target''} we obtain the
connected sum of the lens space $L(p,1)$ and the lens space $L(r,1)$ as
a twisted face-pairing manifold, where $p$ and $r$ are positive
integers.  In Example~\ref{ex:tet1a} we obtain all integer Dehn
surgeries on the figure eight knot as twisted face-pairing manifolds.
In Example~\ref{ex:nil} we obtain the Heisenberg manifold, the prototype
of Nil geometry.  In Example~\ref{ex:s2xs1'} we obtain $S^2\times S^1$.

Which orientable connected closed 3-manifolds
are twisted face-pairing manifolds?
As far as we know they all are, although that seems rather unlikely.  An
interesting problem is to determine whether the 3-torus is a twisted
face-pairing manifold; we do not know whether it is or not.  In a later
paper \cite{survey} we present a survey of twisted face-pairing
3-manifolds which indicates the scope of the set of twisted
face-pairing manifolds.  Here are some of the results in \cite{survey}.
We show how to obtain every lens space as a twisted face-pairing
manifold.  We consider the faceted 3-balls for which every face is a
digon, and we show that the twisted face-pairing manifolds obtained from
these faceted 3-balls are Seifert fibered manifolds.
We show how to obtain most Seifert fibered manifolds.
We show that if $M_1$ and $M_2$ are twisted
face-pairing manifolds, then so is the connected sum of $M_1$ and $M_2$.

This research was supported in part by National Science Foundation grants
DMS-9803868, DMS-9971783, and DMS-10104030. We thank the referee for
helpful suggestions on improving the exposition.

\section{Generalizing the construction}\label{sec:generalizing}\nosubsections

Our twisted face-pairing construction begins with a faceted 3-ball.  In
Section 2 of \cite{twisted} we define a faceted 3-ball $P$ to be an
oriented regular CW complex such that $P$ is a closed 3-ball and $P$ has a
single 3-cell.  In this section we generalize our twisted face-pairing
construction by generalizing the notion of faceted 3-ball. This
generalization gives us more freedom in constructing twisted face-pairing
manifolds, and it is natural in the context of Theorem~\ref{thm:cylinder}.

We take cells of cell complexes to be closed unless explicitly stated
otherwise.

\begin{figure}[ht!]
\centerline{\includegraphics{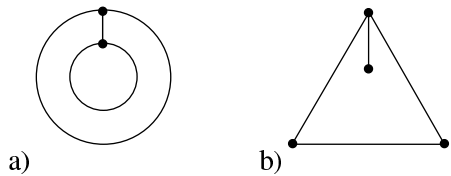}} \caption{Faces now allowed in a
faceted 3-ball} \label{fig:faces}
\end{figure}

We now define a faceted 3-ball $P$ to be an oriented CW complex such that
$P$ is a closed 3-ball, the interior of $P$ is the unique open 3-cell of
$P$, and the cell structure of $\partial P$ does not consist of just one
0-cell and one 2-cell.  Suppose that $P$ is an oriented CW complex such
that $P$ is a closed 3-ball and the interior of $P$ is the unique open
3-cell of $P$.  The condition that the cell structure of $\partial P$ does
not consist of just one 0-cell and one 2-cell is equivalent to the
following useful condition.  For every 2-cell $f$ of $P$ there exists a CW
complex $F$ such that $F$ is a closed disk, the interior of $F$ is the
unique open 2-cell of $F$, and there exists a continuous cellular map
$\zv\co F\to f$ such that the restriction of $\zv$ to every open cell of
$F$ is a homeomorphism.  So $f$ is gotten from $F$ by identifying some
vertices and identifying some pairs of edges.  The number of vertices and
edges in $F$ is uniquely determined.  This definition of faceted 3-ball
allows for faces such as those in Figure~\ref{fig:faces}, which were not
allowed before; part a) of Figure~\ref{fig:faces} shows a quadrilateral
and part b) of Figure~\ref{fig:faces} shows a pentagon.  To overcome
difficulties presented by faces such as those in Figure~\ref{fig:faces},
the next thing that we do is subdivide $P$.

\begin{figure}[ht!]
\centerline{\includegraphics{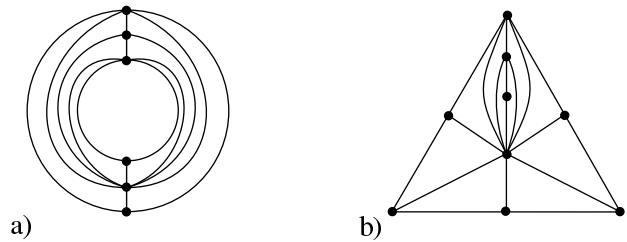}} \caption{The subdivisions of the
faces in Figure~\ref{fig:faces}} \label{fig:subdivision}
\end{figure}

In this paragraph we construct a subdivision $P_s$ of a given faceted
3-ball $P$.  The idea is to not subdivide the 3-cell of $P$ and to
construct what might be called the barycentric subdivision of $\partial
P$.  The vertices of $P_s$ are the vertices of $P$ together with a
barycenter for every edge of $P$ and a barycenter for every face of $P$.
Every face of $P_s$ is a triangle contained in $\partial P$.  If $t$ is
one of these triangles, then one vertex of $t$ is a vertex of $P$, one
vertex of $t$ is a barycenter of an edge of $P$, and one vertex of $t$ is
a barycenter of a face of $P$.  The only 3-cell of $P_s$ is the 3-cell of
$P$.  This determines $P_s$.  Given a face $f$ of $P$, we let $f_s$ denote
the subcomplex of $P_s$ which consists of the cells of $P_s$ contained in
$f$.  Figure~\ref{fig:subdivision} shows $f_s$ for each of the faces $f$
in Figure~\ref{fig:faces}.

In this paragraph we make two related definitions.  Let $P$ be a faceted
3-ball, and let $f$ be a face of $P$.  We define a \textit{corner} of
$f$ at a vertex $v$ of $f$ to be a subcomplex of $f_s$ consisting of the
union of two faces of $f_s$ which both contain an edge $e$ such that $e$
contains $v$ and the barycenter of $f$.  We define an \textit{edge cone}
of $f$ at an edge $e$ of $f$ to be a subcomplex of $f_s$ consisting of
the union of two faces of $f_s$ which both contain an edge $e'$ such
that $e'$ contains the barycenter of $f$ and the barycenter of $e$.

A face-pairing $\ze$ on a given faceted 3-ball $P$ now consists of the
following.  First, the faces of $P$ are paired: for every face $f$ of $P$
there exists a face $f^{-1}\ne f$ of $P$ such that $(f^{-1})^{-1} = f$.
Second, the faces of $P_s$ are paired: for every face $t$ of $P_s$
contained in a face $f$ of $P$ there exists a face $t^{-1}$ of $P_s$ with
$t^{-1}\subseteq f^{-1}$ such that $(t^{-1})^{-1} = t$.  Third, for every
face $t$ of $P_s$ there exists a cellular homeomorphism $\ze_t\co t\to
t^{-1}$ called a partial face-pairing map such that $\ze_{t^{-1}} =
\ze_t^{-1}$.  We require that $\ze_t$ maps the vertex of $P$ in $t$ to the
vertex of $P$ in $t^{-1}$, that $\ze_t$ maps the edge barycenter in $t$ to
the edge barycenter in $t^{-1}$, and that $\ze_t$ maps the face barycenter
in $t$ to the face barycenter in $t^{-1}$. Furthermore, the faces of $P_s$
are paired and the partial face-pairing maps are defined so that if $t$
and $t'$ are faces of $P_s$ contained in some face $f$ of $P$ and if $e$
is an edge of $t\cap t'$ which contains the barycenter of $f$, then
$\ze_t|_e = \ze_{t'}|_e$.  For every face $f$ of $P$ we set $\ze_f =
\{\ze_t\co t\text{ is a face of }f_s\}$, and we refer to $\ze_f$ as a
multivalued face-pairing map from $f$ to $f^{-1}$. We set $\ze =
\{\ze_f\co f\text{ is a face of }P\}$.  In a straightforward way we obtain
a quotient space $P_s/\ze $ consisting of orbits of points of $P _s$ under
$\ze$.  Finally, as in \cite[Section 2]{twisted} we impose on $\ze$ the
face-pairing compatibility condition that as one goes through an edge
cycle the composition of face-pairing maps is the identity.
The cell structure of $P$ induces a cell
structure on $P_s/\ze $, and it is this cell structure that we put on
$P_s/\ze $, not the cell structure induced from $P_s$.  We usually write
$P/\ze $ instead of $P_s/\ze $.  We usually want $\ze$ to be orientation
reversing, which means that every partial face-pairing map of $\ze$
reverses orientation.

Let $P$ be a faceted 3-ball, let $f$ be a face of $P$, and suppose that
$\ze$ is an orientation-reversing face-pairing on $P$.  Then the
multivalued face-pairing map $\ze_f$ determines a function from the set
of corners of $f$ to the set of corners of $f^{-1}$ in a straightforward
way.  The image of one corner of $f$ under this function determines the
image of every corner of $f$ under this function.  The action of $\ze$
on the set of corners of the faces of $P$ determines $P_s/\ze $ up to
homeomorphism.  Thus for our purposes to define the multivalued
face-pairing map $\ze_f$ of a face $f$ of $P$, it suffices to give a
corner $c$ of $f$ and the corner of $f^{-1}$ to which $\ze_f$ maps $c$.

Let $\ze$ be an orientation-reversing face-pairing on a faceted 3-ball
$P$.  Essentially as in Section 2 of \cite{twisted}, $\ze$ partitions the
edges of $P$ into edge cycles.  (We consider the edges of $P$, not the
edges of $P_s$.)  To every edge cycle $E$ of $\ze$ we associate a length
$\ell_E$ and a multiplier $m_E$ as before.  The function $\text{mul}\co
\{\text{edge cycles}\}\to \bN $ defined by $E\mapsto m_E$ is called the
multiplier function.  We obtain a twisted face-pairing subdivision $Q$
from $P$ just as before: if $e$ is an edge of $P$ and if $E$ is the edge
cycle of $\ze$ containing $e$, then we subdivide $e$ into $\ell_Em_E$
subedges.  As before, we subdivide in an $\ze$-invariant way.  We likewise
construct $Q_s$ in an $\ze$-invariant way.  It follows that $\ze$
naturally determines a face-pairing on $Q$, which we continue to call
$\ze$, abusing notation more than before.

We consider face twists in this paragraph.  In the present setting a face
twist is not a single cellular homeomorphism, but instead a collection of
cellular homeomorphisms.  For this, we maintain the situation of the
previous paragraph.  Let $f$ be a face of $Q$.  Let $t$ be a face of
$f_s$.  The orientation of $f$ determines a cyclic order on the faces of
$f_s$.  Let $t'$ be the second face of $f_s$ which follows $t$ relative to
this cyclic order.  Let $\zt_t$ be an orientation-preserving cellular
homeomorphism from $t$ to $t'$ such that $\zt_t$ fixes the barycenter of
$f$.  We call $\zt_f = \{\zt_t\co t\text{ is a face of }f_s\}$ the face
twist of $f$.  We assume that if $t_1$ and $t_2$ are faces of $f_s$ and if
$e$ is an edge of $t_1\cap t_2$ which contains the barycenter of $f$, then
$\zt_{t_1}|_e = \zt_{t_2}|_e$.  We also assume that our face twists are
defined $\ze$-invariantly: for each face $f$ of $Q$ and each face $t$ of
$f_s$, we have $\zt_{t^{-1}} = \ze_{t''}\circ \zt_{t''}^{-1}\circ
\ze_{t^{-1}}$, where $t''$ is the second face of $f_s$ which precedes $t$.
We furthermore impose a compatibility condition on our face twists in the
next paragraph.

Now we are prepared to define a twisted face-pairing $\zd$ on $Q$.  We
pair the faces of $Q$ just as the faces of $P$ are paired.  The pairing of
the faces of $P_s$ likewise induces a pairing of the faces of $Q_s$. For
every face $f$ of $Q$ and every face $t$ of $f_s$, we set $\zd_t =
\ze_{t'}\circ \zt_t$, where $t'$ is the second face of $f_s$ which follows
$t$.  For every face $f$ of $Q$ we set $\zd_f = \{\zd_t\co t\text{ is a
face of }f_s\}$, and we set $\zd = \{\zd_f\co f\text{ is a face of }Q\}$.
We assume that the maps $\zt_t$ are defined so that $\zd$ satisfies the
face-pairing compatibility condition that as one goes through a cycle of
edges in $Q$ the compositions of face-pairing maps is the identity.
Then $\zd$ is a face-pairing on $Q$ called the twisted face-pairing.

Finally, we define $M=M(\ze,\text{mul})$ to be the quotient space
$Q_s/\zd $.  We emphasize that for a cell structure on $M$ we take the
cell structure induced from $Q$, not the cell structure induced from
$Q_s$.  The cell complex $M$ is determined up to homeomorphism by the
function mul and the action of $\ze$ on the corners of the faces of
$P$.

Let $P$ be a faceted 3-ball, let $\ze$ be an orientation-reversing
face-pairing on $P$, and let mul be a multiplier function for $\ze$.
The results of \cite{twisted} all hold in this more general setting.  So
$M$ is an orientable closed 3-dimensional manifold with one vertex.  The
dual of the link of that vertex is isomorphic to $\partial Q^*$ as
oriented 2-complexes, where $Q^*$ is a faceted 3-ball gotten from $Q$ by
reversing orientation.  We label and direct the faces and edges of $Q$
and $Q^*$ as before.  We again obtain a duality between $M$ and $M^*$.
The proofs in \cite{twisted} are valid in the present more general
setting with only straightforward minor technical modifications and the
following.  To obtain a duality between $M$ and $M^*$ in \cite{twisted},
we construct a dual cap subdivision $Q_\zs$ of $Q$.  We let
$C_1,\dotsc,C_k$ be the 3-cells of $Q_\zs$, and for every
$i\in\{1,\dotsc,k\}$ we let $A_i$ be a cell complex isomorphic to $C_i$
so that $A_1,\dotsc,A_k$ are pairwise disjoint.  Then the vertices of
$Q$ can be enumerated as $x_1,\dotsc,x_k$ so that $C_i$ is the unique
3-cell of $Q_\zs$ which contains $x_i$ for $i\in\{1,\dotsc,k\}$.
If $x_i$ has valence $v_i$, then $A_i$ is an alternating suspension on a
$2v_i$-gon for $i\in\{1,\dotsc,k\}$.  In the present setting the
3-cells of $Q_\zs$ need not be alternating suspensions; they are
quotients of alternating suspensions.  See Section~\ref{sec:structure} for
a discussion of the 3-cells of $Q_\zs$.  So in the present
setting we let $x_1,\dotsc,x_k$ be the vertices of $Q$ with valences
$v_1,\dotsc,v_k$, and for $i\in\{1,\dotsc,k\}$ we simply define $A_i$ to
be an alternating suspension on a $2v_i$-gon.  As in \cite{twisted} the
twisted face-pairing $\zd$ on $Q$ induces in a straightforward way what
might be called a face-pairing on the disjoint union of
$A_1,\dotsc,A_k$.  At this point we proceed as in \cite{twisted}.

We conclude this section with
two simple examples to illustrate some of the new phenomena
which occur for our more general faceted 3-balls.

\begin{figure}[ht!]
\centerline{\includegraphics{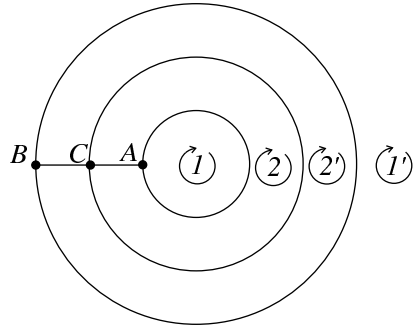}} \caption{The complex $P$ for
Example~\ref{ex:target}} \label{fig:targetp}
\end{figure}

\begin{ex}\label{ex:target} Let the model faceted 3-ball
$P$ be as indicated in Figure~\ref{fig:targetp} with two monogons and
two quadrilaterals, the outer monogon being at infinity.  The inner
monogon has label 1 and is directed outward.  The outer monogon has
label 1 and is directed inward.  The inner quadrilateral has label 2 and
is directed outward.  The outer quadrilateral has label 2 and is
directed inward.  As usual for faces in figures, all four faces are
oriented clockwise.  We construct an orientation-reversing face-pairing
$\ze$ on $P$ as follows.  Multivalued face-pairing map $\ze_1$ maps the
inner monogon to the outer monogon, there being essentially only one way
to do this.  Multivalued face-pairing map $\ze_2$ maps the inner
quadrilateral to the outer quadrilateral fixing their common edge.  Set
$\ze = \{\ze_1^{\pm 1},\ze_2^{\pm 1}\}$.

We might view this face-pairing as follows.  Construct a monogon in the
open northern hemisphere of the 2-sphere $S^2$, put a vertex on the equator
of $S^2$ and join the two vertices with an edge.  Now vertically project
this cellular decomposition of the northern hemisphere into the southern
hemisphere.

The edge cycles for $\ze$ have the following diagrams.
  \begin{equation*}\linnum\label{lin:targetedge}
CC\xrightarrow{\ze_2}CC \qquad
AC\xrightarrow{\ze_2}BC\xrightarrow{\ze_2^{-1}}AC \qquad
BB\xrightarrow{\ze_2^{-1}}AA\xrightarrow{\ze_1}BB
\end{equation*}
For now let the first edge cycle have multiplier 4, let the second have
multiplier 1, and let the third have multiplier 1.

\begin{figure}[ht!]
\centerline{\includegraphics{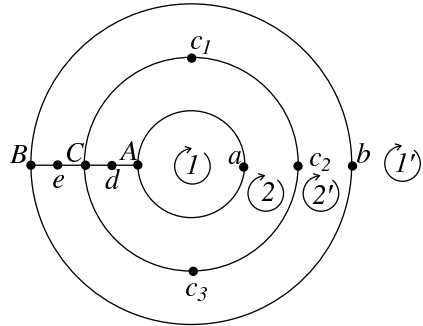}} \caption{The complex $Q$ for
Example~\ref{ex:target}} \label{fig:targetq}
\end{figure}

Figure~\ref{fig:targetq} shows the faceted 3-ball $Q$; we labeled the
new vertices of $Q$ arbitrarily.  Figure~\ref{fig:targetlink} shows
the link of the vertex of $M$, with conventions as in \cite{twisted}.
Figure~\ref{fig:targetqstar} shows the faceted 3-ball $Q^*$ dual to
$Q$ with its edge labels and directions.  Note that $\partial Q^*$ is
dual to the link of the vertex of $M$.  We obtain a presentation for
the fundamental group $G$ of $M$ as follows. Corresponding to the
face labels 1 and 2 we have generators $x_1$ and $x_2$.  The boundary
of the face of $Q^*$ labeled 1 and directed outward gives the relator
$x_1x_2^{-1}$.  The boundary of the face of $Q^*$ labeled 2 and
directed outward gives the relator $x_2^5x_1^{-1}$.  So
  \begin{equation*}
G\cong \langle x_1,x_2\co x_1x_2^{-1},x_2^5x_1^{-1} \rangle \cong \bZ/4
\bZ.
\end{equation*}

\begin{figure}[ht!]
\centerline{\includegraphics{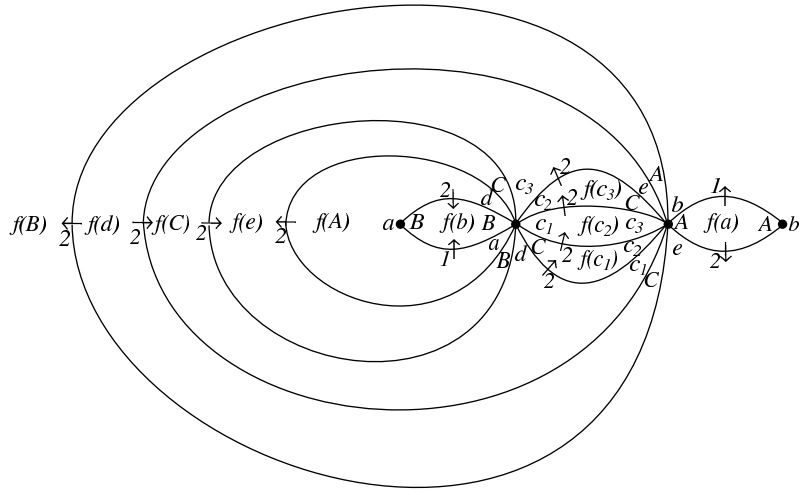}} \caption{The link of the vertex
of $M$} \label{fig:targetlink}
\end{figure}

\begin{figure}[ht!]
\centerline{\includegraphics{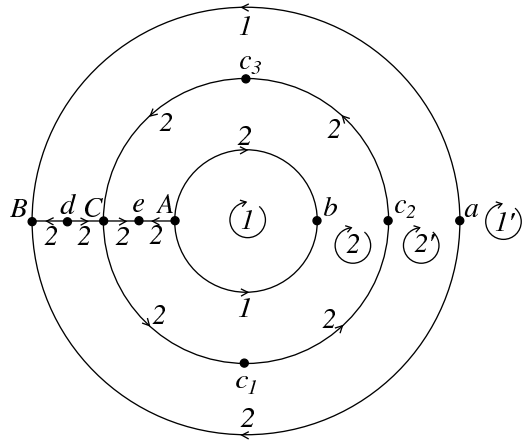}} \caption{The complex $Q^*$
with edge labels and directions} \label{fig:targetqstar}
\end{figure}

We will see in Example~\ref{ex:target''} that $M$ is the lens space
$L(4,1)$.  In general, if the first edge cycle of $\ze$ has multiplier
$p$, if the second edge cycle of $\ze$ has multiplier $q$, and if the
third edge cycle of $\ze$ has multiplier $r$, then we will see in
Example~\ref{ex:target''} that $M$ is the connected sum of the lens space
$L(p,1)$ and the lens space $L(r,1)$ (and so in particular $M$ does not
depend on $q$).
\end{ex}

\begin{figure}[ht!]
\centerline{\includegraphics{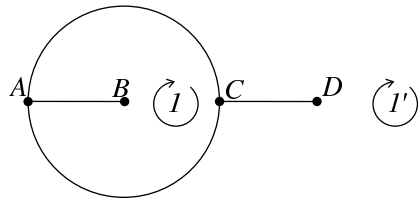}} \caption{The complex $P$
for Example~\ref{ex:s2xs1}}
\label{fig:s2xs1p}
\end{figure}

\begin{ex}[]\label{ex:s2xs1}Let the model faceted 3-ball $P$ be
as in Figure~\ref{fig:s2xs1p} with two quadrilaterals, the outer
quadrilateral being at infinity.  The inner quadrilateral has label 1
and is directed outward.  The outer quadrilateral has label 1 and is
directed inward.  The orientation-reversing multivalued face-pairing map
$\ze_1$ maps the inner quadrilateral to the outer quadrilateral taking
vertex $C$ to vertex $D$.  Set $\ze = \{\ze_1^{\pm 1}\}$.

The vertices $A$ and $C$ of $P$ are joined by two edges.  We use the
subscripts $u$ and $d$ for up and down to distinguish them.  So $AC_u$
is the upper edge joining $A$ and $C$, and $AC_d$ is the lower edge
joining $A$ and $C$.  The face-pairing $\ze$ has only one edge cycle,
and this edge cycle has the following diagram.
  \begin{equation*}
AC_u\xrightarrow{\ze_1}CD\xrightarrow{\ze_1^{-1}}AC_d\xrightarrow{\ze_1^{-1}}
BA\xrightarrow{\ze_1}AC_u
\end{equation*}
For simplicity let this edge cycle have multiplier 1.

Figure~\ref{fig:s2xs1q} shows the faceted 3-ball $Q$; we labeled the
new vertices of $Q$ arbitrarily.  Figure~\ref{fig:s2xs1link} shows
the link of the vertex of $M$.  Figure~\ref{fig:s2xs1qstar} shows the
faceted 3-ball $Q^*$ dual to $Q$ with its edge labels and directions.
Note that $\partial Q^*$ is dual to the link of the vertex of $M$. We
obtain a presentation for the fundamental group $G$ of $M$ as
follows.  Corresponding to the face label 1 we have a generator
$x_1$. The boundary of the face of $Q^*$ labeled 1 and directed
outward gives the trivial relator.  So $G$ has one generator and no
relators, that is, $G\cong \bZ $.

\begin{figure}[ht!]
\centerline{\includegraphics{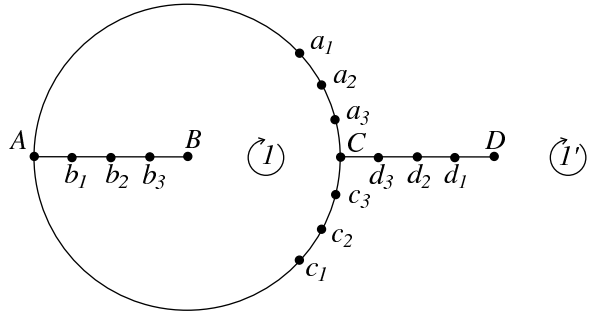}} \caption{The complex $Q$
for Example~\ref{ex:s2xs1}}
\label{fig:s2xs1q}
\end{figure}

\begin{figure}[ht!]
\centerline{\includegraphics{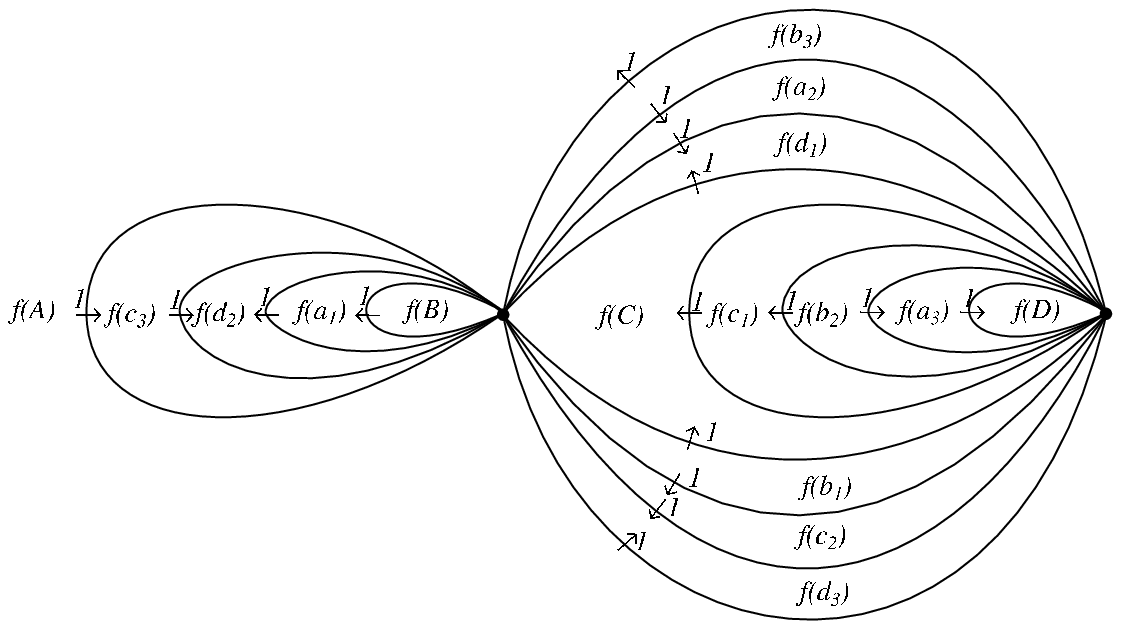}} \caption{The link of the
vertex of $M$} \label{fig:s2xs1link}
\end{figure}

\begin{figure}[ht!]
\centerline{\includegraphics{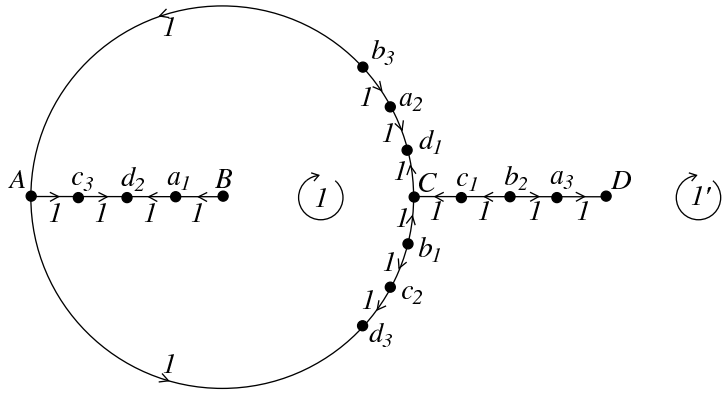}} \caption{The complex $Q^*$
with edge labels and directions} \label{fig:s2xs1qstar}
\end{figure}

We will see in Example~\ref{ex:s2xs1'} that $M$ is homeomorphic to
$S^2\times S^1$ for every choice of multiplier for the edge cycle of
$\ze$.
\end{ex}

\section{Dual cap subdivision}\label{sec:dualcap}\subsections

\subsection{Definition}\label{sec:definition}\nosubsubsections

Recall that we discussed dual cap subdivision in Section 4 of
\cite{twisted}.  Of course, there our faceted 3-balls are regular.  We
generalize to our present cell complexes in a straightforward way.

Let $P$ be a faceted 3-ball.  We construct a dual cap subdivision $P
_\zs$ of $P$ as follows.  The vertices of $P_\zs$ consist of the
vertices of the subdivision $P_s$ defined in
Section~\ref{sec:generalizing} together with a barycenter for the 3-cell
of $P$.  We next describe the edges of $P_\zs$.

The edges of $\partial P_\zs $ consist of the edges of $P_s$ which do
not join the barycenter of a face of $P$ and a vertex of that face.  For
every face of $P$, the subdivision $P_\zs$ also contains an edge joining
the barycenter of that face and the barycenter of the 3-cell of $P$.
These are all the edges of $P_\zs$.

\begin{figure}[ht!]
\centerline{\includegraphics{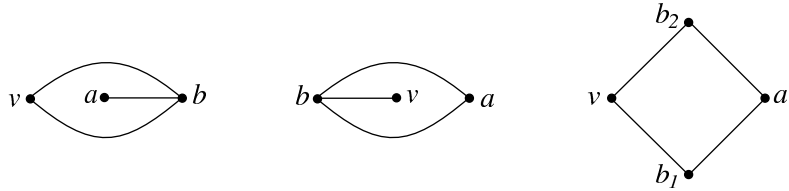}} \caption{The three types of
faces of $\partial P_\zs$} \label{fig:bdyface}
\end{figure}

Having described the edges of $P_\zs$, the structure of $\partial P_\zs
$ is determined.  The faces of $\partial P_\zs $ are in bijective
correspondence with the corners of the faces of $P$.  Every face of
$\partial P_\zs $ is a quadrilateral whose underlying space equals the
underlying space of a corner $c$ at a vertex $v$ of a face $f$ of $P$.
Of course, this quadrilateral contains the barycenter $a$ of $f$.  The
first diagram in Figure~\ref{fig:bdyface} shows this quadrilateral if
$c$ has three vertices and $f$ is a monogon.  The second diagram in
Figure~\ref{fig:bdyface} shows this quadrilateral if $c$ has three
vertices and $f$ is not a monogon.  The third diagram in
Figure~\ref{fig:bdyface} shows this quadrilateral if $c$ has four
vertices. In the first two diagrams $b$ is the barycenter of the edge of
$f$ that contains $v$, and in the third diagram $b_1$ and $b_2$ are the
barycenters of the two edges of $f$ that contain $v$.

\begin{figure}[ht!]
\centerline{\includegraphics{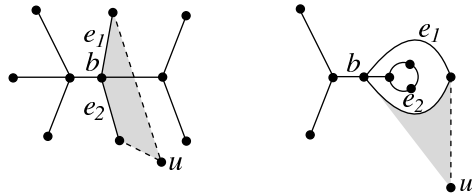}} \caption{Faces of $P_\zs$ not
contained in $\partial P_\zs$} \label{fig:interiorface}
\end{figure}

The remaining faces of $P_\zs$ are in bijective correspondence with the
edges of $P$.  Let $e$ be an edge of $P$, and let $b$ be the barycenter
of $e$.  We have constructed exactly two edges $e_1$ and $e_2$ in
$\partial P_\zs $ which contain $b$ and are not contained in $e$.  The
edge $e$ determines a quadrilateral face of $P _\zs$ containing $e_1\cup
e_2$ and the barycenter $u$ of the 3-cell of $P$.  If $e$ is contained
in two distinct faces of $P$, then the face of $P_\zs$ determined by $e$
has four distinct edges as in the first diagram
of Figure~\ref{fig:interiorface}.
If $e$ is contained in just one face of $P$, then the face of $P_\zs$
determined by $e$ is a degenerate quadrilateral as in the second
diagram of
Figure~\ref{fig:interiorface}.  We have now described all the faces of
$P_\zs$.  This determines $P_\zs$. Note that every vertex of $P$ is in a
unique 3-cell of $P_\zs$.

Now that we have defined dual cap subdivisions of faceted 3-balls, we
define dual cap subdivisions of more general cell complexes.  Let $X$
be a CW complex which is the union of its 3-cells, and suppose that
for every 3-cell $C$ of $X$ there exists a faceted 3-ball $B$ and a
continuous cellular map $\zv\co B\to C$ such that the restriction of
$\zv$ to every open cell of $B$ is a homeomorphism.  We say that a
subdivision $X_\zs$ of $X$ is a dual cap subdivision of $X$ if for
every such choice of $C$ the cell structure on $C$ induced from
$X_\zs$ pulls back via $\zv$ to give a dual cap subdivision of $B$.

It is now clear how to also define a dual cap subdivision of every CW
complex with dimension at most 2 such that every 2-cell contains
an edge. If $X$ is a cell complex for which we have defined a dual cap
subdivision and $k$ is a positive integer, then we let $X_{\zs^k}$
denote the $k$-th dual cap subdivision of $X$.

\subsection{Structure of 3-cells}\label{sec:structure}
\nosubsubsections

In this subsection we discuss the structure of the 3-cells which occur
in the dual cap subdivision of a faceted 3-ball.

\begin{figure}[ht!]
\centerline{\includegraphics{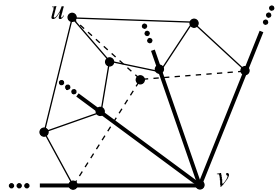}} \caption{The 3-cell $B$ of $P_\zs$
which contains the vertex $v$ of $P$} \label{fig:3cell}
\end{figure}

Let $P$ be a regular faceted 3-ball.  In Section 4 of \cite{twisted} we
showed that every 3-cell of $P_\zs$ is an alternating suspension.  Every
3-cell of $P_\zs$ contains exactly one vertex of $P$, and every vertex
of $P$ is contained in exactly one 3-cell of $P_\zs$.  If $v$ is a
vertex of $P$ with valence $k$, then the 3-cell $B$ of $P_\zs$ which
contains $v$ is an alternating suspension of a $2k$-gon.  See
Figure~\ref{fig:3cell}, which is the same as Figure 15 of
\cite{twisted}.  In Figure~\ref{fig:3cell} the vertex $v$ is a vertex of
$P$ and $u$ is the barycenter of the 3-cell of $P$.
Figure~\ref{fig:3cell} shows an alternating suspension of an octagon.

\begin{figure}[ht!]
\centerline{\includegraphics{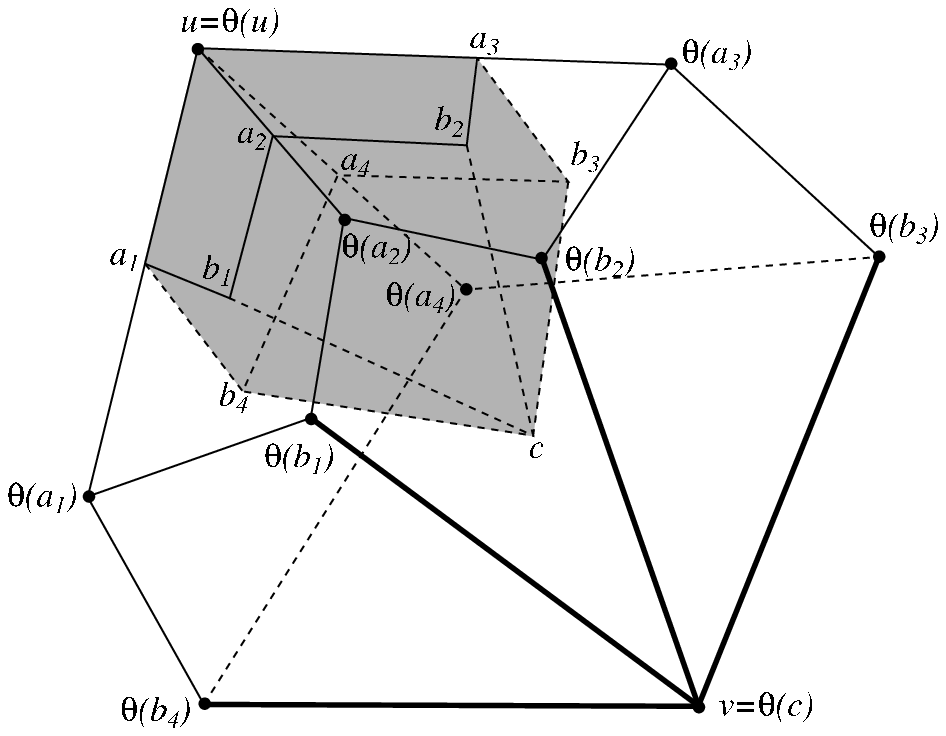}} \caption{
$\text{Star}(u,B_\zs)$} \label{fig:3cellsub}
\end{figure}

We point out here an important property of the dual cap subdivision of an
alternating suspension. Let $B$ be an alternating suspension as in the
previous paragraph. Because the faces of $\partial B$ are quadrilaterals
and $B$ is homeomorphic to the cone on $\text{star}(v,\partial B)$, the
3-cell of $B_\zs$ which contains $u$ is homeomorphic to $B$ by a cellular
homeomorphism $\theta\co \text{star}(u,B_\zs)\to B$ with the following
property: if $x$ is a vertex of $\text{star}(u,B_\zs)$ and $X$ is a cell
of $B$ with $x\in X$, then $\theta(x)\in X$. Figure~\ref{fig:3cellsub}
shows $\text{star}(u,B_\zs)$ for the 3-cell $B$ from Figure~\ref{fig:3cell}.
For convenience further in this section, we make the following definition.
Suppose that $V$ is a CW complex with dimension at most 3,
$U$ is a subcomplex of $V_\zs$, and $\theta\co U\to V$ is a cellular
homeomorphism. We say that $\theta$ \textit{keeps vertices in their cells}
if $\theta(x)\in X$ whenever $x$ is a vertex of $U$ and $X$ a cell in $V$
with $x\in X$.

Now we consider the case of a general faceted 3-ball $P$.  Let $v$ be
a vertex of $P$.  Let $e_1,\dotsc,e_k$ be the edges of $P_\zs$ which
contain $v$.  For every $i\in\{1,\dotsc,k\}$ let $v_i$ be the vertex
of $e_i$ unequal to $v$.  There are $k$ corners of faces at $v$.  Let
$f_1,\dotsc,f_k$ be the faces which contain these corners.  Let $u$
be the barycenter of the 3-cell of $P$, and let $u_i$ be the
barycenter of $f _i$ for every $i\in\{1,\dotsc,k\}$.  If
$u_1,\dotsc,u_k$ and $v_1,\dotsc,v_k$ are distinct, then just as in
the previous paragraph, there is exactly one 3-cell of $P_\zs$ which
contains $v$ and this 3-cell is an alternating suspension of a
$2k$-gon with cone points $u$ and $v$.  In general exactly one 3-cell
of $P_\zs$ contains $v$ and every 3-cell of $P_\zs$ contains exactly
one vertex of $P$.  The 3-cell $C$ of $P_\zs$ which contains $v$ is a
quotient of an alternating suspension $B$ of a $2k$-gon with cone
points mapping to $u$ and $v$, the identifications arising as
follows.  If $f_i=f_j$ for some $i,j\in\{1,\dotsc,k\}$, then
$u_i=u_j$, and so the edge joining $u$ and $u_i$ equals the edge
joining $u$ and $u_j$.  If $v_i=v_j$ for some $i,j\in\{1,\dotsc,k\}$,
then the face containing $u$ and $v_i$ equals the face containing $u$
and $v_j$.  So the 3-cell of $P_\zs$ which contains $v$ is a quotient
of an alternating suspension of a $2k$-gon with cone points mapping
to $u$ and $v$. The quotient map performs two kinds of
identifications.  Edges containing the cone point which maps to $u$
are identified if some face of $P$ is not locally an embedded disk at
$v$.  Faces containing the cone point which maps to $u$ are
identified if some edge of $P$ is not locally an embedded line
segment at $v$.  In every case the restriction of the quotient map to
every open cell of the alternating suspension is a homeomorphism.
Since the identifications are along edges and faces containing $u$,
the map $\theta\co\text{star}(u',B_\zs)\to B$ can be defined so that
it induces a cellular homeomorphism $\zj_C\co \text{star}(u,C_\zs)\to
C$.

\subsection{Central balls}\label{sec:balls}
\nosubsubsections

In this subsection and the next we investigate the second dual cap
subdivision of a faceted 3-ball.

Let $P$ be a faceted 3-ball.  Let $u$ be the vertex of $P_\zs$ which is
the barycenter of the 3-cell of $P$.  Let $C$ be a 3-cell of $P_\zs$.
Section~\ref{sec:structure} shows that $C$ contains $u$,  $C$ is a
quotient of an alternating suspension $B$,
and there is a cellular homeomorphism $\zj_C\co \text{star}(u,C_\zs)\to C$
which keeps vertices in their cells. These homeomorphisms can be
defined compatibly on the pairwise intersections of their domains so that
they piece together to give 
a cellular homeomorphism $\zj\co
\text{star}(u,P_{\zs^2})\to P_\zs $ which keeps vertices in their cells.
We call $\text{star}(u,P_{\zs^2})$ the \textit{central ball} of $P_{\zs^2}$.
We have just shown that the central ball of $P_{\zs^2}$ is cellularly
homeomorphic to $P_\zs$ in a way which is canonical on vertices.

\subsection{Chimneys}\label{sec:chimneys}\nosubsubsections

Let $P$ be a faceted 3-ball.  Let $u$ be the vertex of $P_\zs$ which is
the barycenter of the 3-cell of $P$.  Let $A_1$ be the star of $u$ in
the 1-skeleton of $P_\zs$.  Let $A=\text{star}(A_1,P_{\zs^2})$.  We call
$A$ the \textit{chimney assembly} for $P$.  This subsection is devoted to
investigating the structure of chimney assemblies.

Let $f$ be a face of $P$, and let $a$ be the vertex of $P_\zs$ which is
the barycenter of $f$.  Then $\text{star}(a,P_{\zs^2})$ is a subcomplex
of $A$, which we call the $f$-\textit{chimney} of $A$.

Let $f$ be a face of $P$.  Let $F$ be a CW complex such that $F$ is a
closed disk, the interior of $F$ is the unique open 2-cell of $F$,
and there exists a continuous cellular map $\zv\co F\to f$ such that
the restriction of $\zv$ to every open cell of $F$ is a
homeomorphism. Given a dual cap subdivision $f_\zs$ of $f$, we choose
a dual cap subdivision $F_\zs$ of $F$ so that $\zv$ induces a
cellular map $\zv_\zs\co F_\zs\to f_\zs $.  Let $C_f$ be the mapping
cylinder of $\zv_\zs$, viewed as a CW complex in the obvious way.

\begin{figure}[ht!]
\centerline{\includegraphics{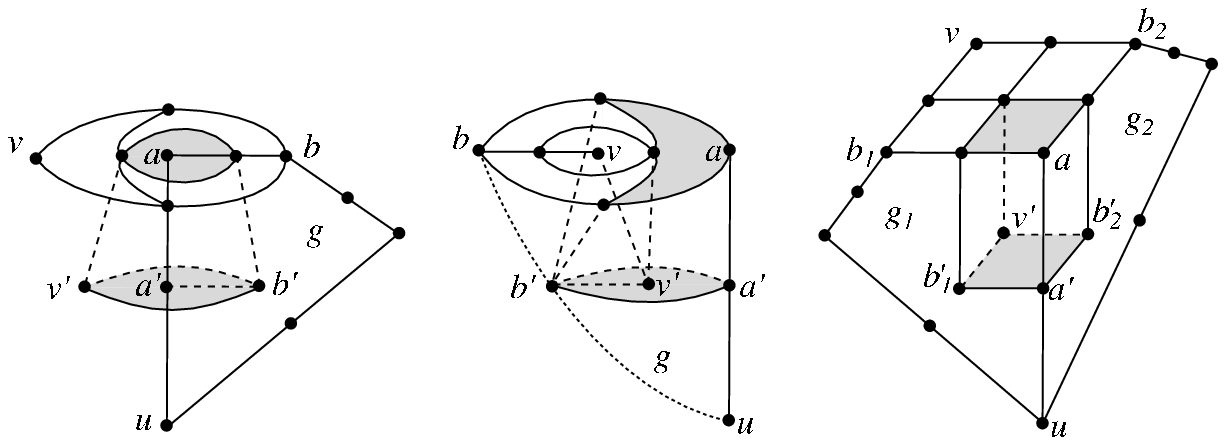}} \caption{Part of $P_{\zs^2}$}
\label{fig:chimneypart}
\end{figure}

In this and the next four
paragraphs we show that $C_f$ is cellularly homeomorphic to the
$f$-chimney of $A$.  Let $a$ be the barycenter of $f$ and let $v$ be a vertex
of $f$. Recall from Figure~\ref{fig:bdyface} and the discussion in
Section~\ref{sec:definition} that there are
three possibilities for a face of $\partial P_{\zs}$.
For each of the three possibilities, Figure~\ref{fig:chimneypart} shows part of
$P_{\zs^2}$.  Every vertex and edge in Figure~\ref{fig:chimneypart} is a
vertex or edge of $P_{\zs^2}$ except for the dotted arc in the second
diagram which joins $b$, $b'$, and $u$. The barycenter $a$ of $f$ is shown.
In the first two diagrams $b$ is the barycenter of the edge of $f$
that contains $v$, and $a$, $b$, and $v$ are the vertices
of a face $h$ of $f_\zs$. In the third diagram $b_1$ and $b_2$ are
the barycenters of the two edges of $f$ that contain $v$, and
$a$, $b_1$, $b_2$, and $v$ are the vertices of a face $h$ of $f_\zs$.
The dual cap subdivision of $h$ is shown in Figure~\ref{fig:chimneypart}.
The barycenter $u$ of $P$ and $a$ are
joined by an edge $e$ of $P_\zs$.  Let $a'$ be the barycenter of $e$ in
$P_{\zs^2}$.  Let the map $\zj\co \text{star}(u,P_{\zs^2})\to P_\zs $ be
as in Section~\ref{sec:balls}.  Section~\ref{sec:balls} shows that
$\zj(a')=a$. Let $C$ be the 3-cell of
$P_\zs$ which contains $v$, and let $v'$ be the barycenter of $C$ in
$P_{\zs^2}$. Section~\ref{sec:balls} shows that $\zj(v')=v$.
Let $k$ be the face of $h_\zs$ which contains $a$.
In each of the three diagrams in Figure~\ref{fig:chimneypart}
we have drawn in gray the face $k$ and a face $h'$ which will
be described below. We consider separately
the three possibilities for $h$ shown in Figure~\ref{fig:bdyface}.

We first consider the case that $h$ has the form of the first diagram
in Figure~\ref{fig:bdyface}. Then $f$ is a monogon.
Let $g$ be the face of $P_\zs$ which
contains $a$, $b$, and $u$, and let $b'$ be the barycenter of $g$.
For clarity, two edges of $g_{\zs}$ are not shown.
Section~\ref{sec:balls} shows that $\zj(b')=b$. Let $h'$ be the face of
$P_{\zs^2}$ with vertices $a'$, $b'$, and $v'$.
Then $k$ and $h'$ are cellularly homeomorphic,
$\text{star}(a,P_{\zs^2})$ is the product of a 1-simplex and the dual cap
subdivision of a monogon, and $\text{star}(a,P_{\zs^2})$ is cellularly
homeomorphic to $C_f$.

Now suppose that $h$ has the form of the second diagram
in Figure~\ref{fig:bdyface}. Then $v$ has valence $1$ in $\partial f$.
As in the previous case let $g$ be the face of $P_\zs$ which
contains $a$, $b$, and $u$, and let $b'$ be the barycenter of $g$.
Section~\ref{sec:balls} again shows that $\zj(b')=b$. Let $h'$ be the face of
$P_{\zs^2}$ with vertices $a'$, $b'$, and $v'$. Then
$k$ is cellularly homeomorphic to a square and $h'$ is cellularly
homeomorphic to a square with two adjacent edges identified.  It follows
that the 3-cell of $\text{star}(a,P_ {\zs^2})$ which contains $v'$ is
cellularly homeomorphic to a cube with two adjacent edges identified.

Finally, suppose $h$ has the form of the third diagram in Figure
\ref{fig:bdyface}. For $i\in\{1,2\}$, let $g_i$ be the face of
$P_\zs$ which contains $u$ and $b_i$, and let $b'_i$ be the vertex of
$P_{\zs^2}$ which is the barycenter of $g_i$. For clarity two edges
of $(g_1)_\zs $ and two edges of $(g_2)_\zs $ are omitted in the
third diagram in Figure~\ref{fig:chimneypart}.
 Section~\ref{sec:balls} shows that
$\zj(b'_i)=b_i$ for $i\in\{1,2\}$.
Let $h'$ be the face of $P_{\zs^2}$ with vertices
$a'$, $b_1'$, $b_2'$ and $v'$. We see that $\zj$ restricts to a cellular
homeomorphism from $h'$ to $h$.
Then both
$k$ and $h'$ are cellularly homeomorphic to squares and the 3-cell of
$\text{star}(a,P_{\zs^2})$ which contains $v'$ is cellularly homeomorphic
to a cube.

If $h$ has the form of the second or third diagram in
Figure~\ref{fig:bdyface}, then $\text{star}(a,P_{\zs^2})$ is a union of
complexes as described in the previous two paragraphs.
It follows in these cases that $\text{star}(a,\partial
P_{\zs^2})$ is cellularly homeomorphic to $F_\zs$, that the restriction
of $\zj$ to $\text{star}(a,P_{\zs^2})\cap\text{star}(u,P_{\zs^2})$ is a
cellular homeomorphism onto $f_\zs$, and that $\text{star}(a,P_{\zs^2})$ is
cellularly homeomorphic to $C_f$.

\begin{figure}[ht!]
\centerline{\includegraphics{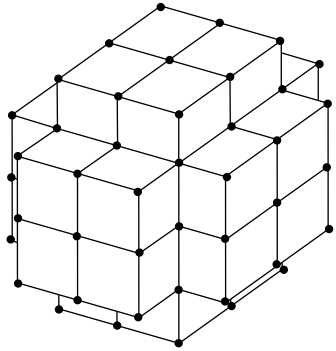}} \caption{The chimney assembly
for a cube} \label{fig:dualcapstar}
\end{figure}

So the chimney assembly $A$ for $P$ is the union of the central ball of
$P_{\zs^2}$ and the chimneys of the faces of $P$.  The central ball of
$P_{\zs^2}$ is cellularly homeomorphic to $P_\zs$, and the chimneys of
the faces of $P$ are mapping cylinders.  Figure~\ref{fig:dualcapstar}
shows the chimney assembly for a cube.

Let $f$ be a face of $P$, and let $C_f$ be the $f$-chimney of $A$.  We
call $f\cap C_f$ the \textit{top} of $C_f$.  We call the intersection of
$C_f$ with the central ball of $A$ the \textit{bottom} of $C_f$.  We
call faces of $\partial C_f$ which are in neither the top nor the bottom
of $C_f$ \textit{lateral faces}.

\section{Building Heegaard diagrams from face-pairings}
\label{sec:heegaard}\subsections

In this section we construct Heegaard diagrams from face-pairings.

\subsection{Edge pairing surfaces}\label{sec:edge}\nosubsubsections

We begin by constructing a cellulated closed surface $S$ from a
face-pairing.  We call $S$ the \textit{edge pairing surface} of the
face-pairing.  See the introduction, where $S$
is defined for regular faceted 3-balls.  Our more general faceted
3-balls present some complications, but we proceed in much the same way.

Let $P$ be a faceted 3-ball with orientation-reversing face-pairing $\ze$.
We first construct a cell complex $X$ cellularly homeomorphic to the
1-skeleton of $P$.  Let $f$ and $f^{-1}$ be two paired faces of $P$. Next
construct a CW complex $F$ such that $F$ is a closed disk, the interior of
$F$ is the unique open 2-cell of $F$, and there exists a continuous
cellular map $\zv\co F\to f$ such that the restriction of $\zv$ to every
open cell of $F$ is a homeomorphism.  There also exists a corresponding
cellular map $\zj\co F\to f^{-1}$ such that $\zv$ and $\zj$ are related as
follows.  Recall that to define $\ze$ we construct subdivisions $f_s$ and
$f_s^{-1}$ of $f$ and $f^{-1}$ in Section~\ref{sec:generalizing}.  Let $t$
be a face of $f_s$.  Then there exists a corresponding face $t^{-1}$ of
$f_s^{-1}$ and a partial face-pairing map $\ze_t\co t\to t^{-1}$.  There
also exists a subspace $T$ of $F$ such that the restriction of $\zv$ to
$T$ is a homeomorphism onto $t$.  We may, and do, choose the maps $\zv$
and $\zj$ so that if $x\in T$, then $\zj(x)=\ze_t(\zv(x))$.  We next
construct $\partial F\times[0,1]$.  We view the interval $[0,1]$ as a
1-cell, and we view $\partial F\times[0,1]$ as a 2-complex with the
product cell structure. For every $x\in\partial F$ we identify
$(x,0)\in\partial F\times[0,1]$ with the point of $X$ corresponding to
$\zv(x)\in\partial f$ and we identify $(x,1)\in\partial F\times[0,1]$ with
the point of $X$ corresponding to $\zj(x)\in\partial f^{-1}$.  Doing this
for every pair of faces of $P$ yields a cell complex $Y$ whose underlying
space is a closed surface.  We define $S$ to be the dual cap subdivision
of $Y$. We say that an edge of $S$ is \textit{vertical} if it is either
contained in $X$ or is disjoint from $X$.  We say that an edge of $S$ is
\textit{diagonal} if it is not vertical.  We say that an edge of $S$ is a
\textit{meridian edge} if it is not an edge of $Y$.  We refer to edges of
$Y$ as \textit{nonmeridian} edges of $S$.
The union of the vertical meridian edges is a family
$\{\za_1,\dotsc,\za_n\}$ of pairwise disjoint simple closed curves in $S$
called the \textit{vertical meridian curves} of $S$.
The union of the diagonal meridian edges is a family
$\{\zb_1,\dotsc,\zb_m\}$ of pairwise disjoint simple closed curves in $S$
called the \textit{diagonal meridian curves} of $S$.

\begin{figure}[ht!]
\centerline{\includegraphics{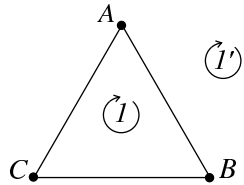}} \caption{The complex $P$
for Example~\ref{ex:lens}}
\label{fig:triangle}
\end{figure}

\begin{figure}[ht!]
\centerline{\includegraphics{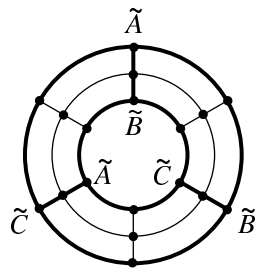}} \caption{The edge pairing surface
of $\ze$ viewed as a quotient of an annulus} \label{fig:lenss}
\end{figure}

\begin{figure}[ht!]
\centerline{\includegraphics{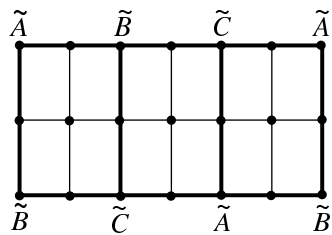}} \caption{The edge pairing
surface of $\ze$ viewed as a quotient of a quadrilateral}
\label{fig:threequadrilaterals}
\end{figure}

\begin{figure}[ht!]
\centerline{\includegraphics{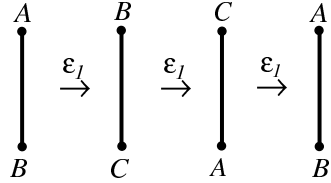}} \caption{A diagram of the edge
cycle of $\ze$} \label{fig:edgecycle}
\end{figure}

\begin{ex}[]\label{ex:lens} We illustrate the above edge pairing
surface construction using the simple example of the lens space
$L(3,1)$.  To obtain $L(3,1)$ we take a faceted 3-ball $P$ with just two
faces which are triangles as in Figure~\ref{fig:triangle}, where one
face is at infinity.  The orientation-reversing face-pairing map $\ze_1$
maps the inner triangle to the outer triangle taking vertex $A$ to
vertex $B$.  We set $\ze = \{\ze_1^{\pm 1}\}$.  Let $S$ be the edge
pairing surface of $\ze$, and let $\widetilde{A}$, $\widetilde{B}$, and
$\widetilde{C}$ be the vertices of $S$ which correspond to $A$, $B$, and
$C$.  Figure~\ref{fig:lenss} shows $S$ as an annulus whose boundary
components are to be identified in a straightforward way.  Similarly,
Figure~\ref{fig:threequadrilaterals} shows $S$ as a quotient of a
quadrilateral.  This quadrilateral is gotten from the edge cycle of
$\ze$, shown in Figure~\ref{fig:edgecycle}, in a straightforward way.
The meridian edges of $S$ are drawn with thin arcs, and the
nonmeridian edges of $S$ are drawn with thick arcs.  We see
that $S$ is a torus.  The union of the vertical meridian edges is a
simple closed curve on the torus, and the union of the diagonal meridian
edges is a simple closed curve on the torus.  The torus and these two
curves form a Heegaard diagram for $L(3,1)$.  This is a special case of
Theorem~\ref{thm:heegaard}.
\end{ex}

\subsection{Heegaard diagrams for general face-pairings}
\label{sec:general}\nosubsubsections

\begin{thm}[]\label{thm:heegaard}
Let $P$ be a faceted 3-ball with orientation-reversing face-pairing $\ze$
such that the cell complex $N=P/\ze$ is a manifold with one vertex.
Let $H_1$ be the star of the barycenter of the 3-cell of $N$ in the
1-skeleton of $N_\zs$, and let $H=\text{star}(H_1,N_{\zs^2})$.  Then
\begin{description}
\item[i)]$H$ is a handlebody in $N$
and $\partial H$ is a Heegaard surface for $N$, and
\item[ii)]
$\partial H$ is cellularly homeomorphic to the edge pairing surface $S$
of $\ze$.
\end{description}
Identifying $S$ with $\partial H$, we have the following:
\begin{description}
\item[iii)] the set $\{\za_i\}$ of vertical meridian curves of $S$
is a basis of meridian curves for $H$;
\item[iv)] the set $\{\zb_i\}$ of diagonal meridian curves of $S$
is a basis of meridian curves for $N\setminus \text{int}(H)$;
\item[v)] $(S; \{\za_i\};\{\zb_i\})$ is a Heegaard diagram for $N$.
\end{description}
\end{thm}
\begin{proof}
We view $N_{\zs^2}$ as a quotient of $P_{\zs^2}$.  The preimage
of $H$ in $P_{\zs^2}$ is the chimney assembly $A$ for $P$, and $H$ is
obtained from the chimney assembly by gluing together the tops in pairs.
Hence $H$ is a handlebody in $N$. Let $H'$ denote the closure of the
complement of $H$ in $N$. Then $H'$ is the star in $N_{\zs^2}$ of the
1-skeleton of $N$, and so is a union of stars of vertices and stars of
edge barycenters. Each vertex star is a cone with cone point the vertex,
and hence is homeomorphic to a closed ball since $N$ is a manifold. Each
3-cell in the star of an edge barycenter is either a cube or (as in the
second diagram in Figure~\ref{fig:chimneypart}) a cube with a pair of adjacent
edges identified, and the star of each barycenter is a two-sided mapping
cylinder obtained from $D^2 \times I$. It follows that $H'$ is a
handlebody and so $\partial H$ is a Heegaard surface for
$N$. This proves i).

In this paragraph we show that $\partial H$ is cellularly homeomorphic
to $S$.   The
preimage of $\partial H$ in $P_{\zs^2}$ is the union of all the lateral
faces of the chimneys of $A$.  Section~\ref{sec:chimneys} shows that
every chimney of $A$ is a mapping cylinder, and so the union of the
lateral faces of every chimney of $A$ is a mapping cylinder.
Hence $\partial H$ is homeomorphic to a topological space obtained by
attaching two-sided mapping cylinders to the 1-skeleton of $P$. It now
follows from the definition of $S$ in terms of attaching two-sided mapping
cylinders that $S$ is cellularly homeomorphic to $\partial H$. This proves
ii).

Let $f$ be a face of $P$.  The top of the $f$-chimney $C_f$ of $A$ meets
the union of the lateral faces of $C_f$ in a simple closed edge path in
$A$.  This edge path maps to a meridian curve for $H$, and every edge in
this meridian curve corresponds to a vertical meridian edge of $S$.
Hence the union of the edges of $\partial H$ corresponding
to the vertical meridian edges of $S$ forms a basis of meridian curves
for $H$. This proves iii).

Suppose given an edge cycle of $\ze$ consisting of $j$ distinct edges
$e_1,\dotsc,e_j$ of $P$ with diagram
  \begin{equation*}
e_1\xrightarrow{\ze_{f_1}}e_2\xrightarrow{\ze_{f_2}}\cdots
\xrightarrow{\ze_{f_{j-1}}}e_j\xrightarrow{\ze_{f_j}}e_1.
\end{equation*}
Let $u$ be the vertex of $P_{\zs^2}$ which is the barycenter of the 3-cell
of $P$,  and let $\zj\co \text{star}(u,P_{\zs^2})\to P_\zs $ be the cellular
homeomorphism of Section~\ref{sec:balls}.  Let $e'_i=\zj^{-1}((e_i)_\zs)$,
let $v_i$ be the vertex of $e'_i$ such that $\zj(v_i)$ is the barycenter
of $e_i$ and let $C_{f_i}$ be the $f_i$-chimney of $A$ for every
$i\in\{1,\dotsc,j\}$.  For every $i\in\{1 ,\dotsc,j\}$ the chimney
$C_{f_i}$ contains two lateral faces and the chimney $C_{f_i^{-1}}$
contains two lateral faces with the following properties, where $i+1$
is taken modulo $j$.  See Figure~\ref{fig:chimneysides}.
The two lateral faces of $C_{f_i}$ both
contain an edge which contains $v_i$ and a vertex $x _i$ in the top of
$C_{f_i}$, and the two lateral faces of $C_{f_i^{-1}}$ both contain an
edge which contains $v_{i+1}$ and a vertex $y_i$ in the top of
$C_{f_i^{-1}}$.  Furthermore the image in $\partial H$ of $x_i$ equals the
image in $\partial H$ of $y_i$, and both the edge containing $v_i$ and
$x_i$ and the edge containing $v_{i+1}$ and $y_i$ map to edges of
$\partial H$ which correspond to
diagonal meridian edges of $S$.
For each $i\in\{1,\dotsc,j\}$, $y_{i-1}$, $v_i$, $x_i$, and the barycenter
$b_i$ of $e_i$ are vertices of face of a
$P_{\zs^2}$, where $i-1$ is taken modulo
$j$. The union in $N_{\zs^2}$ of the images of these $j$ faces is
a properly embedded closed disk $D$ in $H'$ whose boundary is in
$\partial H' = \partial H$. See Figure~\ref{fig:mercurve}, which shows
a chimney assembly together with the top of an $f$-chimney drawn with
thick arcs and the face of $P_{\zs^2}$ with vertices $y_{i-1}$, $v_i$,
$x_i$, and $b_i$.
If $\ze$ has $m$ edge cycles, then we
obtain $m$ such disks $D_1,\dotsc,D_m$ in $H'$.
The disks $D_1,\dotsc,D_m$ are pairwise disjoint,
$H'\setminus\bigcup_{i=1}^{m}D_i$
has one connected component for every vertex of $N$, and each of these
connected components contracts to the corresponding vertex.  Since $N$ has
only one vertex, it follows that $D_1,\dotsc,D_m$ form a basis of meridian
disks for $H'$, and so the union of the edges of $\partial H$
corresponding to the diagonal meridian edges of $S$ forms a basis of
meridian curves for $H'$. This proves iv), and v) now follows.

\begin{figure}[ht!]
\centerline{\includegraphics{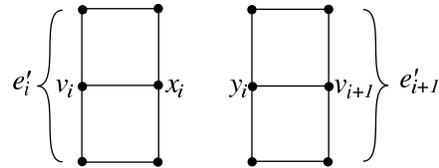}} \caption{Two lateral faces of
$C_{f_i}$ and two lateral faces of $C_{f^{-1}_i}$}
\label{fig:chimneysides}
\end{figure}

\begin{figure}[ht!]
\centerline{\includegraphics{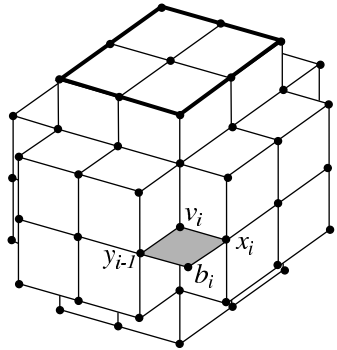}} \caption{
The face with vertices $y_{i-1}$, $v_i$, $x_i$, and $b_i$}
\label{fig:mercurve}
\end{figure}

This proves Theorem~\ref{thm:heegaard}.
\end{proof}

\subsection{Heegaard diagrams for twisted face-pairing 3-manifolds}
\label{sec:twisted}\nosubsubsections

We next interpret Theorem~\ref{thm:heegaard} for twisted
face-pairing 3-manifolds.

Let $P$ be a faceted 3-ball with orientation-reversing face-pairing $\ze$,
and suppose given a multiplier function for $\ze$.  Let $Q$ be the
associated twisted face-pairing subdivision of $P$, let $\zd$ be the
associated twisted face-pairing on $Q$, and let $M=Q/\zd $ be the
associated twisted face-pairing manifold.  Let $S$ be the edge pairing
surface of $\zd$.

\begin{figure}[ht!]
\centerline{\includegraphics{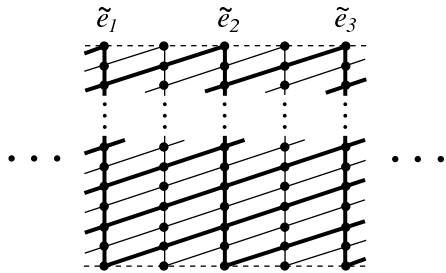}} \caption{The edge cycle
cylinder corresponding to the $\ze$-edge cycle $E$}
\label{fig:edgecylinder}
\end{figure}

Theorem~\ref{thm:heegaard} implies that $S$ is cellularly homeomorphic
to a Heegaard surface for $M$.  We view $S$ as a union of subspaces, one
for every edge cycle of $\ze$ as follows.  Let $E$ be an edge cycle of
$\ze$.  Suppose that $E$ has length $j$, multiplier $k$ and edge cycle
diagram
  \begin{equation*}
e_1\xrightarrow{\ze_{f_1}}e_2\xrightarrow{\ze_{f_2}}\cdots
\xrightarrow{\ze_{f_{j-1}}}e_j\xrightarrow{\ze_{f_j}}e_1.
\end{equation*}
To construct $Q$ from $P$ we subdivide each of the edges $e_1,\dotsc,e_j$
into $jk$ subedges.  Every edge of $Q$ gives rise to two edges of $S$.  So
the edges $e_1,\dotsc,e_j$ of $P$ give rise to subcomplexes
$\widetilde{e}_1,\dotsc,\widetilde{e}_j$ of $S$ each of which is the union
of $2jk$ edges of $S$.  As in Figure 11 of \cite{twisted}, $\zd$ maps
subedge $m$ of $e_i$ relative to $f_i$ to subedge $m+1$ of $e_{i+1}$
relative to $f_{i+1}$ for every $i\in\{1,\dotsc,j\}$ and
$m\in\{1,\dotsc,jk-1\}$, where $i+1$ is taken modulo $j$.  It follows that
$E$ gives rise to a subspace $C$ of $S$ as shown in
Figure~\ref{fig:edgecylinder}.  We call $C$ an \textit{edge cycle
cylinder}.  Certain arcs contained in $C$ are not edges of $S$, and so
they are drawn with dashes.  The edges of $S$ are drawn with two
thicknesses simply to distinguish the thin meridian edges from the thick
nonmeridian edges of $S$.  In general $C$ need not be homeomorphic to a
closed annulus, but there exists a closed annulus $A$ and a surjective
continuous map $\zv\co A\to C$ such that the restriction of $\zv$ to the
interior of $A$ is a homeomorphism and $\zv$ maps the boundary of $A$ to
the union of the arcs drawn with dashes in Figure~\ref{fig:edgecylinder}.
We refer to the images under $\zv$ of the two boundary components of $A$
as the \textit{ends} of $C$.  The ends of $C$ are chosen so that the edge
cycle cylinders corresponding to different $\ze$-edge cycles meet only
along their boundaries and their union is $S$.  If $\zg$ is an arc in $A$
which joins the boundary components of $A$, then we say that the curve
$\zv(\zg)$ \textit{joins the ends} of $C$.  If $\zg$ is a simple closed
curve in the interior of $A$ which separates the boundary components of
$A$, then we say that the curve $\zv(\zg)$ \textit{separates the ends} of
$C$.  We define the \textit{circumference} of $C$ to be $j$, and we define
the \textit{height} of $C$ to be $jk$.  Now we see that Figure 11 of
\cite{twisted} essentially shows an edge cycle cylinder in a Heegaard
surface for $M$.  The thick vertical edges in
Figure~\ref{fig:edgecylinder} arise from $P$, and the thick diagonal edges
in Figure~\ref{fig:edgecylinder} arise from $P^*$.  Vertical edges of $S$
are drawn vertically, and diagonal edges of $S$ are drawn diagonally.  The
following theorem is now clear.

\begin{thm}[]\label{thm:twistedheegaard}
Let $M = M(\ze,\text{mul})$  be a twisted face-pairing manifold, and
let $\zd$ be the associated twisted face-pairing. Let $S$ be the edge
pairing surface of $\zd$, let $\{\za_i\}$ be the set of vertical
meridian curves, and let $\{\zb_i\}$ be the set of diagonal meridian
curves. Then $(S,\{\za_i\},\{\zb_i\})$ is a Heegaard diagram for $M$.
\end{thm}

\begin{figure}[ht!]
\centerline{\includegraphics{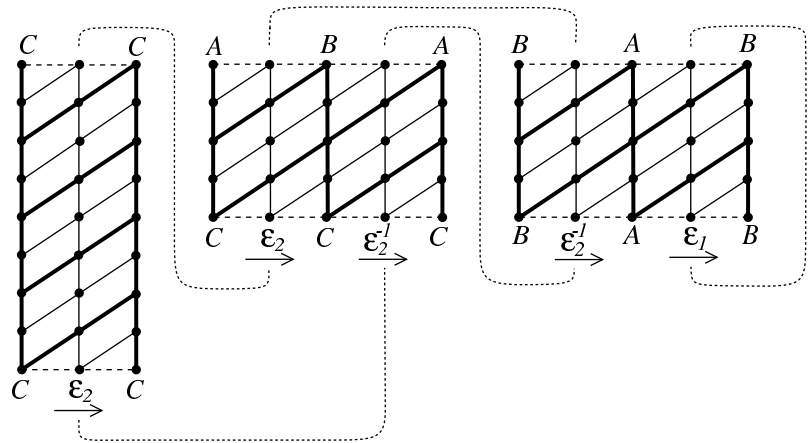}} \caption{A Heegaard diagram
decomposed into three edge cycle cylinders} \label{fig:targetheegaard}
\end{figure}

\begin{ex}[]\label{ex:target'} We return to
Example~\ref{ex:target}.  The model face-pairing in
Example~\ref{ex:target} has three edge cycles.
Line~\ref{lin:targetedge} gives diagrams for them.  As in
Example~\ref{ex:target}, we choose multipliers to be 4, 1, and 1.  Each
of these three edge cycles gives rise to an edge cycle cylinder as in
Figure~\ref{fig:edgecylinder}.  These three edge cycle cylinders are
shown in Figure~\ref{fig:targetheegaard}.  They are drawn as
quadrilaterals with their left sides to be identified with their right
sides.  The first edge cycle cylinder has circumference 1 and height 4,
the second has circumference 2 and height 2, and the third has
circumference 2 and height 2.  The thin dotted arcs in
Figure~\ref{fig:targetheegaard} indicate how the ends of the cylinders
are to be identified.  These identifications respect the face-pairing
maps, which are also shown.  After performing the required
identifications we obtain a closed orientable surface $S$ of genus 2.
The union of its vertical meridian edges is a basis of meridian curves
for $S$, and the union of its diagonal meridian edges is a basis of
meridian curves for $S$.  The result is a Heegaard diagram for our
twisted face-pairing manifold.
\end{ex}

\section{Building face-pairings from Heegaard diagrams}
\label{sec:face-pair}\subsections

In Section~\ref{sec:heegaard} we construct Heegaard diagrams from
face-pairings.  Theorem~\ref{thm:twistedheegaard} shows that every
twisted face-pairing manifold has a Heegaard diagram which can be
decomposed into cylinders which correspond to the edge cycles of the
model face-pairing.  The height of every such cylinder is a multiple of
its circumference, the multiple being the multiplier of the
corresponding edge cycle.  In this section we show that the
decomposition of Heegaard diagrams into analogous cylinders is a general
phenomenon, not one restricted to twisted face-pairing manifolds.  In
general the heights of the cylinders need not be multiples of their
circumferences.  In fact, Theorem~\ref{thm:cylinder} shows that the
height of every such cylinder coming from a given Heegaard diagram is a
multiple of its circumference if and only if the Heegaard diagram arises
from a twisted face-pairing manifold as in
Theorem~\ref{thm:twistedheegaard}.  This provides a characterization of
the Heegaard diagrams which we construct for twisted face-pairing
manifolds.

\subsection{Generalities concerning Heegaard diagrams}
\label{sec:generalities}\nosubsubsections

For us a Heegaard diagram is a Heegaard diagram for a closed
connected orientable
3-manifold.  It consists of an orientable
connected closed surface $S$ with positive
genus and two bases of meridian curves for $S$.  We assume that there
exists a triangulation of $S$ for which each of these meridian curves is
piecewise linear, and we assume that these curves intersect transversely
in only finitely many points. Let $U$ be the union of the two bases of
meridian curves for $S$.  We say that our Heegaard diagram is
\textit{irreducible} if every connected component of $S\setminus U$ is
homeomorphic to an open disk.

Suppose given an irreducible Heegaard diagram consisting of an
orientable connected closed surface $S$ and two bases of meridian curves for $S$.
We refer to the meridian curves in one basis as \textit{vertical
meridian curves}, and we refer to the meridian curves in the other basis
as \textit{diagonal meridian curves}.  The assumptions imply that the
meridian curves of our Heegaard diagram determine a cell structure on
$S$ whose vertices are the intersections of the meridian curves and
whose faces are the closures of the connected components of the
complement in $S$ of the union of the meridian curves.  We refer to the
edges of $S$ which are contained in vertical meridian curves as
\textit{vertical (meridian) edges}, and we refer to the edges of $S$
which are contained in diagonal meridian curves as \textit{diagonal
(meridian) edges}.  Since the meridian curves intersect
transversely, every vertex of $S$ has
valence 4 and the edges of every face of $S$
are alternately vertical and diagonal. Since the Heegaard diagram is
irreducible, no face can have a single edge and so every face of $S$ has an
even number of edges.

\subsection{Heegaard cylinders}\label{sec:cylinders}
\nosubsubsections

Suppose given an irreducible Heegaard diagram with surface $S$.  We view
$S$ as having a cell structure as in the last paragraph.  This
subsection is devoted to defining subspaces of $S$ called Heegaard
cylinders.

\begin{figure}[ht!]
\centerline{\includegraphics{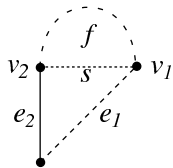}} \caption{The temporary
horizontal segment $s$ of $f$} \label{fig:segment}
\end{figure}

\begin{figure}[ht!]
\centerline{\includegraphics{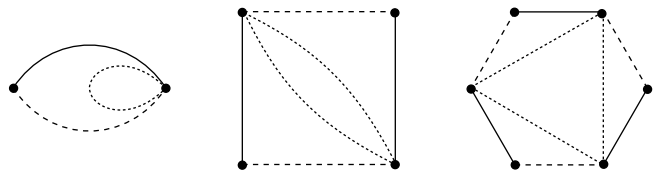}} \caption{A complete set of
temporary horizontal segments for a digon, a quadrilateral and a hexagon}
\label{fig:segments}
\end{figure}

In this paragraph we construct what we call temporary horizontal segments
of $S$.  For this we choose an orientation of $S$.  This orientation of
$S$ determines an orientation of the boundary of every face of $S$.  Let
$f$ be a face of $S$.  Let $v_1$ be a vertex of $f$ such that a diagonal
edge $e_1$ of $f$ follows $v_1$ (relative to $f$).
See Figure~\ref{fig:segment}, where, as usual,
faces are oriented in the clockwise direction. The vertex $v_1$ and
the edge $e_1$ determine a vertical edge $e_2$ of $f$ which follows $e_1$
(relative to $f$) and a terminal vertex $v_2$ of $e_2$ (relative to $f$).
We choose an open arc in the interior of $f$ whose closure joins $v_1$ and
$v_2$.  We call the closure of this open arc a \textit{temporary
horizontal segment} of $S$.  In Figure~\ref{fig:segment},
$e_1$ is drawn with a
dashed line segment, $e_2$ is drawn with a line segment, the rest of the
boundary of $f$ is drawn with a broken arc, and the temporary horizontal
segment $s$ joining $v_1$ and $v_2$ is drawn with a dotted line segment.
We choose a temporary horizontal segment for every such choice of $e_1$
and $e_2$ so that the temporary horizontal segments associated to distinct
choices of $e_1$ and $e_2$ meet only at vertices of $S$.
Figure~\ref{fig:segments} shows a complete set of temporary horizontal
segments for a digon, a quadrilateral, and a hexagon, with conventions as
in Figure~\ref{fig:segment}.

\begin{figure}[ht!]
\centerline{\includegraphics{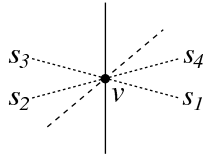}} \caption{A neighborhood of a
vertex $v$ of $S$} \label{fig:neighborhood}
\end{figure}

In this paragraph we define what it means for one temporary horizontal
segment to follow another.  Every vertex $v$ of $S$ has a neighborhood
as in Figure~\ref{fig:neighborhood}.  The vertex $v$ is contained in
temporary horizontal segments $s_1$, $s_2$, $s_3$, and $s_4$, which need not
be distinct.  Rotating about $v$ in the clockwise direction from $s_1$,
we encounter a vertical edge, then a diagonal edge and then $s_2$.  We
say that $s_2$ \textit{follows} $s_1$ and likewise that $s_4$ follows
$s_3$.  If faces are oriented in the counterclockwise direction, then we
rotate about $v$ in the counterclockwise direction.  For every temporary
horizontal segment $s_1$ of $S$ there exists a unique temporary
horizontal segment $s_2$ of $S$ such that $s_2$ follows $s_1$.
Furthermore, $s_1$ is the unique temporary horizontal segment of $S$
such that $s_2$ follows $s_1$.

\begin{figure}[ht!]
\centerline{\includegraphics{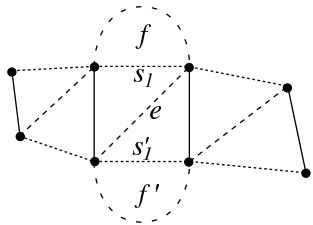}} \caption{The temporary
horizontal segment $s'_1$ of $f'$} \label{fig:mirrorsegment}
\end{figure}

In this paragraph we use the temporary horizontal segments of $S$ to
construct annuli in $S$.  For this let $s_1$ be a temporary horizontal
segment of $S$.  The previous paragraph implies that there exist
temporary horizontal segments $s_2,\dotsc,s_k$ such that $s_{i+1}$
follows $s_i$ for every $i\in\{1,\dotsc,k\}$, where $i+1$ is taken
modulo $k$.  The union of $s_1,\dotsc,s_k$ is a closed curve $\zs$ which
intersects itself at most tangentially, not transversely.  The temporary
horizontal segment $s_1$ is contained in a face $f$ of $S$, and $s_1$ is
related to a diagonal edge $e$ of $f$ as in
Figure~\ref{fig:mirrorsegment}.  Across $e$ from $f$ is a face $f'$ of
$S$, and just as $e$ is related to $s_1$, the edge $e$ is related to a
temporary horizontal segment $s'_1$ in $f'$ as in
Figure~\ref{fig:mirrorsegment}.  Just as $s_1$ determines the closed
curve $\zs$, the temporary horizontal segment $s'_1$ determines a closed
curve $\zs'$.  The curves $\zs$ and $\zs'$ are the boundary
components of an open annulus in $S$ which contains the interior of $e$.

\begin{figure}[ht!]
\centerline{\includegraphics{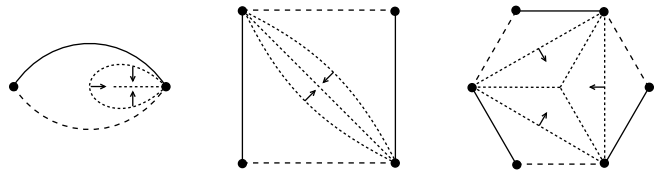}} \caption{Homotoping the
temporary horizontal segments in Figure~\ref{fig:segments}}
\label{fig:homotop}
\end{figure}

A defect of the annuli constructed in the previous paragraph is that the
union of their closures is not all of $S$.  To remedy this defect, we
homotop the temporary horizontal segments of $S$ as indicated in
Figure~\ref{fig:homotop}.  More precisely, for every face $f$ of $S$
choose a barycenter $b$
in the open subset of $f$ bounded by temporary horizontal segments and
join $b$ with an arc to the initial vertex(s)
(relative to $f$) of every diagonal edge of $f$ so that these arcs meet
only at $b$ and they meet the temporary horizontal segments only at
vertices. Then homotop (isotop except for a digon) the temporary
horizontal segments of $S$ contained in $f$ to the union of these arcs,
fixing endpoints.  We refer to the image of a temporary horizontal
segment under such a homotopy as a \textit{horizontal segment}.  The
result of these homotopies is to enlarge the annuli of the previous
paragraph so that the union of their closures is $S$.  We refer to the
closures of these enlarged annuli as \textit{simple cylinders}. Each
simple cylinder $C$ is the image of a closed annulus $A$ under a
continuous map that restricts to a homeomorphism from $\text{int}(A)$ to
$\text{int}(C)$. We call the image of each component of $\partial A$ an
\textit{end} of the simple cylinder. Each end of a simple cylinder is a
union of horizontal segments.

Suppose that $C_1,\dotsc,C_k$ are simple cylinders, and suppose that
$C_i$ has ends $E_i$ and $E'_i$ for every $i\in\{1,\dotsc,k\}$.  Also
suppose that the horizontal segments in $E'_i$ equal the horizontal
segments in $E_{i+1}$ for every $i\in\{1,\dotsc,k-1\}$.  Then we call
$C_1\cup \dots \cup C_k$ a \textit{cylinder}.  We define a \textit{Heegaard
cylinder} to be a cylinder which is maximal with respect to containment.
We define the \textit{height} of a Heegaard cylinder to be the number of
simple cylinders contained in it.  We define the \textit{circumference}
of a Heegaard cylinder to be the number of diagonal edges in any simple
cylinder contained in the given Heegaard cylinder.  The interiors of the
simple cylinders of $S$ are pairwise disjoint, and the union of the simple
cylinders of $S$ is $S$.  It follows that the interiors of the
Heegaard cylinders of $S$ are pairwise disjoint, and the union of the Heegaard
cylinders of $S$ is $S$.

\subsection{Face-pairings for general Heegaard diagrams}
\label{sec:face-pairings}\nosubsubsections

\begin{thm}[]\label{thm:cylinder} Suppose given an irreducible
Heegaard diagram $D$.  Then there exists a faceted 3-ball $P$ with
orientation-reversing face-pairing $\ze$ such that $N=P/\ze $ is a
manifold with one vertex and $D$ is the Heegaard diagram of $N$
described in Theorem~\ref{thm:heegaard}.  Furthermore, $D$ is the
Heegaard diagram of a twisted face-pairing manifold as described in
Theorem~\ref{thm:twistedheegaard} if and only if the height of every
Heegaard cylinder of $D$ is a multiple of its circumference.
\end{thm}
\begin{proof} Let $S$ be the surface of the Heegaard diagram $D$.  We
begin by defining a 1-complex $K$, which is a subspace of $S$.  Recall
that homotoping the temporary horizontal segments to the horizontal
segments in Section~\ref{sec:cylinders} involves choosing a barycenter
for every face of $S$.  These barycenters are the vertices of $K$.  The
edges of $K$ are dual to the diagonal edges of $S$.  In other words, for
every diagonal edge $e$ of $S$ there are faces $f_1$ and $f_2$ of $S$ on
either side of $e$, and there is an edge of $K$ corresponding to $e$
which joins the barycenters of $f_1$ and $f_2$.

Let $V$ be the union of the vertical meridian curves of $D$.  Then
$S\setminus V$ is homeomorphic to the 2-sphere with $2g$ holes, where
$g$ is the genus of $S$.  Of course, we construct $K$ so that
$K\subseteq S\setminus V$. Figure~\ref{fig:fibers}
indicates how to define a strong deformation retraction from
$S\setminus V$ to $K$.  In
Figure~\ref{fig:fibers} horizontal segments are drawn with dotted arcs,
diagonal edges are drawn with dashed arcs, edges of $K$ are drawn with
thick arcs, vertical edges are drawn with medium thick arcs and
retraction fibers are drawn with thin arcs and dashed arcs.
Since filling in the $2g$ holes of $S\setminus V$ yields a 2-sphere,
$K$ is cellularly homeomorphic to the 1-skeleton of a
faceted 3-ball $P$ with $2g$ faces
such that a neighborhood of $K$ in $\partial P$
is homeomorphic to $S\setminus V$.  We identify $K$ with the 1-skeleton
of $P$.

\begin{figure}[ht!]
\centerline{\includegraphics{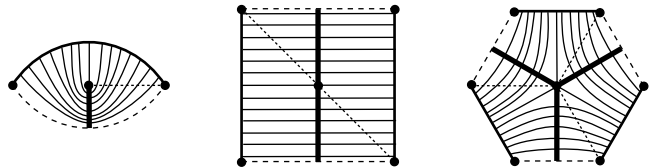}} \caption{Fibers of a retraction
from $S\setminus V$ to $K$} \label{fig:fibers}
\end{figure}

There exists an orientation-reversing face-pairing $\ze$ on $P$ which
acts on the vertices and edges of $K$ as follows.  Let $e$ be an edge of
$K$ with vertices $v_1$ and $v_2$.  By definition $e$ is dual to a
diagonal edge $d$ of $S$.  Let $v$ be a vertex of this diagonal edge of
$S$.  See Figure~\ref{fig:nearv}, where conventions are as in
Figure~\ref{fig:fibers}.  A vertical meridian curve of $D$ passes
through $v$.  Let $d'$ be the diagonal edge of $S$ incident to $v$
across this vertical meridian curve from $d$.  Let $e'$ be the edge of
$K$ dual to $d'$, and let $v'_1$ and $v'_2$ be the
vertices of $e'$ corresponding
to the vertices $v_1$ and $v_2$ of $e$ as in Figure~\ref{fig:nearv}.
The vertex $v$ determines a face $f$ of $P$ which contains $e$
and a face $f^{-1}$ which contains $e'$.  Then
the (multivalued) face-pairing map $\ze_f$ maps $e$ to $e'$ taking $v_1$
to $v'_1$ and $v _2$ to $v'_2$.

\begin{figure}[ht!]
\centerline{\includegraphics{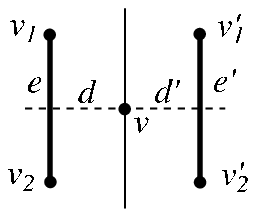}} \caption{Part of $S$ near $v$}
\label{fig:nearv}
\end{figure}

In this paragraph we show that the cell complex $N=P/\ze $ is a manifold
with one vertex.  As in the proof of the main theorem of \cite{intro},
to prove that $N$ is a manifold, it suffices to prove that the Euler
characteristic of $N$ is 0.  It is clear that $N$ has one 3-cell and $g$
faces.  So as in the proof of the main theorem of \cite{intro}, to prove
that $N$ is a manifold, it suffices to prove that $N$ has one vertex and
$g$ edges.  The description of $\ze$ in the previous paragraph shows
that the $\ze$-edge cycles are in bijective correspondence with the
diagonal meridian curves of $D$.  Since $D$ has $g$ diagonal meridian
curves, it follows that $N$ has $g$ edges.  Just as we defined the
1-complex $K$ with edges dual to the diagonal edges of $S$, it is
possible to define a 1-complex $K^*$ with edges dual to the vertical
edges of $S$.  See Figure~\ref{fig:nearv'}, which is the same as
Figure~\ref{fig:nearv}, except that two edges of $K^*$ are added as
thick dashed line segments.  Just as the complex $K$ is connected, so is
the complex $K^*$.  The connectivity of $K^*$ and the description of
$\ze$ in the previous paragraph imply that $N$ has one vertex.  Thus $N$
is a manifold with one vertex.

\begin{figure}[ht!]
\centerline{\includegraphics{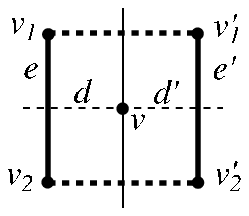}} \caption{Part of $S$ near $v$}
\label{fig:nearv'}
\end{figure}

Figure~\ref{fig:nearv'} indicates how the edges of $K$ and $K^*$
decompose $S$ into a union of quadrilaterals. Furthermore, the
vertical and diagonal median curves of $S$ are edges of the dual cap
subdivision of this tiling by quadrilaterals. Since the two edges of
$K$ that are in a quadrilateral are paired under $\ze$ and $K$ is the
1-skeleton of $P$, this exactly matches the cell structure of the
edge pairing surface $S'$ of $\ze$ as the dual cap subdivision of a
tiling by quadrilaterals. Hence there exists a homeomorphism from $S$
to $S'$ such that the vertical edges of $S$ map to vertical meridian
edges of $S'$ and diagonal edges of $S$ map to diagonal meridian
edges of $S'$.  It follows that $D$ is the Heegaard diagram of $N$
described in Theorem~\ref{thm:heegaard}.

Now suppose that $\ze$ is a twisted face-pairing. By the previous
paragraph, we can identify $S$ with the edge pairing surface $S'$ of
$\ze$. Each edge cycle cylinder is a union of simple cylinders, so
every Heegaard cylinder of $D$ is a union of edge cycle cylinders.
Since the height of every edge cycle cylinder is a multiple of its
circumference, the height of every Heegaard cylinder is a multiple of
its circumference.

Finally suppose that the height of every Heegaard cylinder of $D$ is a
multiple of its circumference. Let $K'$ be the 1-complex with the same
underlying space as $K$ but such that each open edge of $K'$ is a
component of the intersection of $K$ with the interior of a Heegaard
cylinder. Let $P'$ be the faceted 3-ball with the same underlying space as
$P$, but with edges the edges of $K'$ rather than the edges of $K$. Before
we constructed a face-pairing $\ze$ on $P$ by identifying edges of $K$
along diagonal meridian edges. Now we construct a face-pairing $\ze'$ on
$P'$ by identifying edges of $K'$ along horizontal segments. Then two
edges of $K'$ are in the same edge cycle exactly if they are in the same
Heegaard cylinder, and $\ze$ is the twisted face-pairing obtained from
$\ze'$ by choosing the multiplier of an edge cycle to be the quotient of
the height of the corresponding Heegaard cylinder by its circumference.

This proves Theorem~\ref{thm:cylinder}.
\end{proof}

\section{Surgery descriptions for twisted face-pairing manifolds}
\label{sec:surgery}\subsections

The Heegaard diagrams of twisted face-pairing manifolds described in
Theorem~\ref{thm:twistedheegaard} easily yield surgery descriptions for
these manifolds.  This section deals with these surgery descriptions.

\subsection{Initial surgery descriptions}\label{sec:initial}\nosubsubsections

\begin{figure}[ht!]
\centerline{\includegraphics{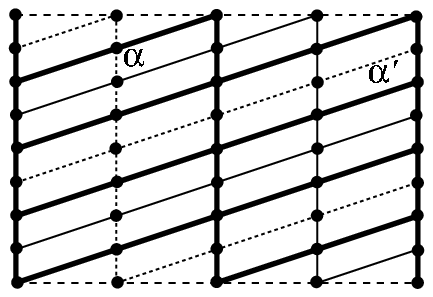}} \caption{The curves $\za$ and
$\za'$} \label{fig:curves}
\end{figure}

Let $P$ be a faceted 3-ball with orientation-reversing face-pairing
$\ze$, and suppose given a multiplier function for $\ze$.  Let $M$ be
the associated twisted face-pairing manifold.
Theorem~\ref{thm:twistedheegaard} describes a Heegaard diagram $D$ for
$M$.  Let $S$ be the surface of $D$, let $C$ be an edge cycle cylinder
of $D$, and let $\za$ be a minimal union of vertical meridian edges of $C$
which joins the ends of $C$.  Figure~\ref{fig:curves} shows $C$ as a
quadrilateral whose left and right sides are to be identified, and $\za$
is shown as a union of vertical dotted edges.  Let $\za'$ be the minimal
union of diagonal meridian edges of $C$ which joins the endpoints of
$\za$ as in Figure~\ref{fig:curves}.  Let $\zb$ be a simple closed curve
in $C$ which separates the ends of $C$.  If the height of $C$ equals the
circumference of $C$, then $\za'$ is isotopic (relative endpoints) to a
Dehn twist of $\za$ along $\zb$.
In this case, let $\zt$ be the appropriate Dehn
twist so that $\za'$ is isotopic (relative endpoints) to $\zt(\za)$.
In general, if the $\ze$-edge cycle corresponding to $C$ has multiplier
$m$, then the height of $C$ divided by the circumference of $C$ equals
$m$ and $\za'$ is isotopic (relative endpoints) to $\zt^m(\za)$
for the appropriate Dehn twist $\zt$ along $\zb$.  In
Figure~\ref{fig:curves}, $m=2$.  This discussion essentially proves the
following theorem.

\begin{thm}[]\label{thm:dehntwist} Let $M = M(\ze,\text{mul})$ be a
twisted face-pairing manifold. Let $S$ be the associated Heegaard surface,
and let $\za_1,\dotsc,\za_n$ be the vertical meridian curves. Let
$\tau^{\text{mul}} = \tau_1^{\text{mul}(E_1)}\circ \dots \circ
\tau_m^{\text{mul}(E_m)}$, where $\tau_1,\dotsc,\tau_m$ are consistently
oriented Dehn twists along the core curves of the edge cycle cylinders
corresponding to the edge cycles $E_1,\dotsc,E_m$. Then
$(S;\{\za_i\},\{\tau^{\text{mul}}(\za_i)\})$ is a Heegaard diagram for $M$.
\end{thm}
\begin{proof} The theorem follows from the previous discussion except
for the matter of the directions of the Dehn twists.  The previous
discussion shows that every edge cycle cylinder determines a Dehn twist.
Since the twisting directions of these Dehn twists are
consistent relative to a fixed orientation of $S$, there are two
choices for $\zt^{\text{mul}}$.  For one choice of $\zt^{\text{mul}}$
the curves $\zt^{\text{mul}}(\za_1),\dotsc,\zt^{\text{mul}}(\za_n)$ are
isotopic to the diagonal meridian curves of $D$, and
Theorem~\ref{thm:dehntwist} is clear.  For the other choice of
$\zt^{\text{mul}}$ the curves
$\zt^{\text{mul}}(\za_1),\dotsc,\zt^{\text{mul}}(\za_n)$ are isotopic to
the diagonal meridian curves, not of $D$, but of the corresponding
Heegaard diagram for the twisted face-pairing manifold $M^*$ dual to
$M$.  Theorem~\ref{thm:dehntwist} follows because Theorem 4.6 of
\cite{twisted} (together with its generalization in
Section~\ref{sec:generalizing} if $P$ isn't regular)
shows that $M^*$ is homeomorphic to $M$.
\end{proof}

We are now prepared for the following theorem, which shows how to obtain
twisted face-pairing manifolds by Dehn surgery on framed links in $S^3$.
We begin with some terminology.

Let $(H,H';\zg_1,\dotsc,\zg_n;\zg_1',\dotsc,\zg_n')$
be a Heegaard splitting of $S^3$. That is, $H$ and $H'$ are handlebodies
in $S^3$ such that $H \cup H' = S^3$ and $H\cap H' = \partial H = \partial
H'$, $\zg_1,\dotsc,\zg_n$ bound a basis of meridian disks for $H$, and
$\zg_1',\dotsc,\zg_n'$ bound a basis of meridian disks for $H'$. We call
$(H,H';\zg_1,\dotsc,\zg_n;\zg_1',\dotsc,\zg_n')$
a \textit{standard Heegaard splitting} of $S^3$ if $\zg_i$ and
$\zg_j'$ are disjoint if $i\ne j$ and $\zg_i \cap \zg_j'$ intersect in a
single point if $i=j$.

If $H$ is a handlebody and $\zb_1,\dotsc,\zb_m$ are pairwise disjoint
simple closed curves in $\partial H$, then pairwise disjoint simple closed
curves $\zd_1,\dotsc,\zd_m$ in the interior of $H$ are \textit{parallel
copies of} $\zt_1,\dotsc,\zt_m$ if they satisfy the following:
there are pairwise disjoint closed annuli $A_1,\dots,A_m$ in $H$
such that $\partial A_i = \zd_i \cup \zt_i$ and
$\zb_i=A_i\cap\partial H$ for every $i\in\{1,\dotsc,m\}$.

\begin{thm}[]\label{thm:surgery}
Let $M = M(\ze,\text{mul})$ be a twisted face-pairing manifold, let
$E_1,\dotsc,E_m$ be the edge cycles of $\ze$,
let $S$ be the edge pairing surface for the twisted face pairing $\zd$,
let $\za_1,\dotsc,\za_n$ be the vertical meridian curves of $S$, and let
$\zb_1,\dotsc,\zb_m$ be core curves for the edge cycle cylinders
corresponding to $E_1,\dotsc,E_m$.
Let $(H,H';\zg_1,\dotsc,\zg_n;\zg_1',\dotsc,\zg_n')$ be a standard
Heegaard diagram of $S^3$.
Let $\zv\co S\to \partial H$ be a homeomorphism such that
$\zv(\za_i)=\zg_i$ for every $i\in\{1,\dotsc,n\}$, and
let $\zd_1,\dotsc,\zd_m$ be parallel copies of
$\zv(\zb_1),\dotsc,\zv(\zb_m)$ in the interior of $H$.
We obtain a link $L$ in $S^3$ by taking
$L=\zg_1\cup\ldots\cup \zg_n\cup \zd_1\cup\ldots\cup \zd_m$.  We define a
framing of $L$ as follows.  The components $\zg_1,\dotsc,\zg_n$ have
framing 0.  For every $i\in\{1,\dotsc,m\}$ the component $\zd_i$ has
framing $\text{lk}(\zd_i,\zv(\zb_i))\pm \text{mul}(E_i)^{-1}$, where
$\text{lk}(\zd_i,\zv(\zb_i))$ is the linking number of $\zd_i$ and
$\zv(\zb_i)$ after they are compatibly oriented and the sign is either $+$
for every $i\in\{1,\dotsc,m\}$ or $-$ for every $i\in\{1,\dotsc,m\}$.
Then the manifold obtained by Dehn surgery on this framed link $L$ is
homeomorphic to $M$.
\end{thm}
\begin{proof} The surface $\partial H$ and the curves
$\zg_1,\dotsc,\zg_n$ and $\zg'_1,\dotsc,\zg'_n$ form a Heegaard diagram
for $S^3$.  By performing Dehn surgery on $\zg_1,\dotsc,\zg_n$, each
with framing 0, we obtain a connected sum of $n$ copies of $S^2\times
S^1$, which has a Heegaard diagram consisting of the surface $\partial
H$, the curves $\zg_1,\dotsc,\zg_n$ and the curves $\zg_1,\dotsc,\zg_n$.
(The bases of meridian curves are equal.)  For every
$i\in\{1,\dotsc,m\}$ let $\zt_i$ be a Dehn twist on $\partial H$ along
$\zv(\zb_i)$, choosing $\zt_1,\dotsc,\zt_m$ so that the directions in
which they twist are consistent relative to a fixed orientation of
$\partial H$.  Set $\zt^{\text{mul}}= \zt_1^{\text{mul}(E_1)}\circ\cdots
\circ\zt_m^{\text{mul}(E_m)}$.  Theorem~\ref{thm:dehntwist}
implies that $M$ has a Heegaard diagram consisting of the surface
$\partial H$, the curves $\zg_1,\dotsc,\zg_n$ and the curves
$\zt^{\text{mul}}(\zg_1),\dotsc,\zt^{\text{mul}}(\zg_n)$.  The fact that
$M$ is obtained by Dehn surgery on $\zg_1,\dotsc,\zg_n$ and
$\zd_1,\dotsc,\zd_m$ now follows from a standard argument which
appears, for example, in the proof of the Dehn-Lickorish Theorem on page
84 of \cite{PS}.  It only remains to determine the framings of
$\zd_1,\dotsc,\zd_m$.

We determine the framings of $\zd_1,\dotsc,\zd_m$ in this paragraph.
Let $i\in\{1,\dotsc,m\}$.  Let $T$ be a solid torus regular neighborhood
of $\zd_i$ such that $\zv(\zb_i)\subseteq\partial T$.
Let $\za\subseteq\partial T$ be
the boundary of a meridian disk of $T$.  The curve $\za$ and part of
$\zv(\zb_i)$ are shown in part a) of Figure~\ref{fig:framing}.  Using
our usual orientation convention for figures as in
Figure~\ref{fig:curves}, our Dehn twist takes $\za$ to a curve $\zg$ as
shown in part b) of Figure~\ref{fig:framing}.  Let $m=\text{mul}(E_i)$.
Then $\zg$ is homologous in $\partial T$ to
$\za-m \zv(\zb_i)= \left(1- m\text{lk}(\zd_i,\zv(\zb_i))\right)\alpha
-m\ell_i$,
where $\ell_i = \zv(\zb_i) - \text{lk}(\zd_i,\zv(\zb_i))\alpha$
is parallel to $\zd_i$ and hence is a longitude for $T$.  Thus the framing of
$\zd_i$ is $\text{lk}(\zd_i,\zv(\zb_i))-1/m$.  If our Dehn twist is in
the opposite direction, then the framing of $\zd_i$ is
$\text{lk}(\zd_i,\zv(\zb_i))+1/m$.

\begin{figure}[ht!]
\centerline{\includegraphics{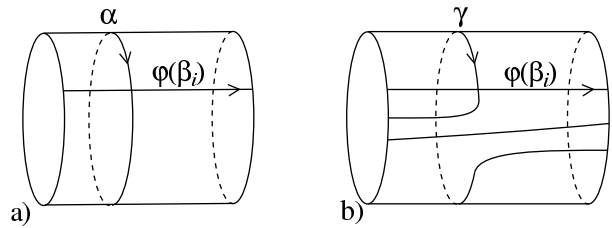}} \caption{Determining the
framing of $\zd_i$} \label{fig:framing}
\end{figure}

This proves Theorem~\ref{thm:surgery}.
\end{proof}

\subsection{The corridor construction}\label{sec:corridor}
\nosubsubsections

Theorem~\ref{thm:surgery} describes a framed link $L$ in $S^3$ such that
Dehn surgery on $L$ obtains a given twisted face-pairing manifold.  The
goal of this subsection is to make the construction of such links $L$
algorithmic and simple.  We call the method which we use the corridor
construction.

\begin{figure}[ht!]
\centerline{\includegraphics{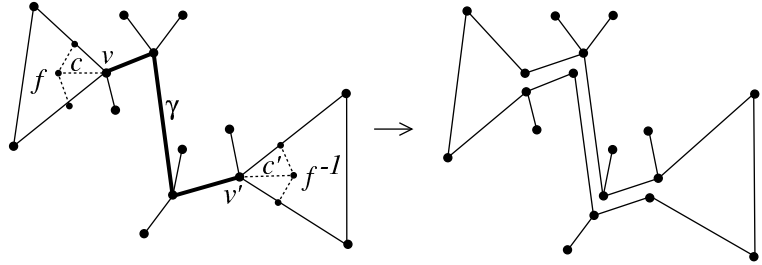}} \caption{Constructing a
corridor between $f$ and $f^{-1}$} \label{fig:corridor}
\end{figure}

Let $P$ be a faceted 3-ball, and let $\ze$ be an
orientation-reversing face-pairing on $P$.  In this paragraph we
construct corridors between the paired faces of $P$.  Let $f$ be a
face of $P$.  The face $f$ is paired with the face $f^{-1}$.  Let $c$
be a corner of $f$ at a vertex $v$ of $f$, and suppose that $\ze_f$
takes $c$ to the corner $c'$ of $f^{-1}$ at the vertex $v'$ of
$f^{-1}$.  Let $\zg$ be an edge path arc in $P$ with endpoints $v$
and $v'$.  See the left part of Figure~\ref{fig:corridor}, where $f$
and $f^{-1}$ are triangles, $\zg$ is drawn with thick line segments,
and the corners $c$ and $c'$ are indicated with dotted edges.  From
$\partial P$ we construct a new cell complex with underlying space
the 2-sphere as follows.  We choose an arbitrarily small neighborhood
of $\zg$ in $\partial P$ and modify the cell structure of $\partial
P$ only in this neighborhood as indicated in
Figure~\ref{fig:corridor}.  The right part of
Figure~\ref{fig:corridor} shows the new cell complex.  We refer to
this modification of $\partial P$ as \textit{constructing a corridor}
between $f$ and $f^{-1}$.  In a straightforward way we continue to
successively construct corridors between all the paired faces of $P$.
We call the resulting cell complex $C$ a \textit{corridor complex}
for $\ze$.  Every face of $C$ is in some sense the union of two
paired faces of $P$ and a corridor.

\begin{figure}[ht!]
\centerline{\includegraphics{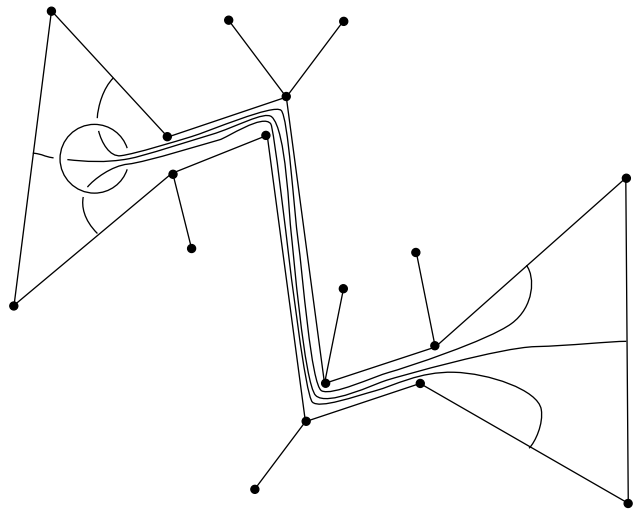}} \caption{The face
component $\za$ and the projections $\zb$} \label{fig:corarcs}
\end{figure}

Again let $P$ be a faceted 3-ball, and let $\ze$ be an
orientation-reversing face-pairing on $P$.  In this paragraph we describe
a planar diagram $D$ of a link $L$ in $S^3$.  Let $C$ be a corridor
complex for $\ze$.  We view the underlying space of $C$ as the one-point
compactification $\bR^2\cup \{\infty\}$ of $\bR^2$, where the point
$\infty$ lies in the interior of some face of $C$.  The diagram $D$ lies
in $C\setminus\{\infty\}$.  Let $g$ be a face of $C$.  We next describe
the part of $D$ which lies in $g$.  One component of $L$ has a projection
$\za$ in the interior of $g\setminus\{\infty\}$ with no self-crossings; it
is unknotted.  We call this component of $L$ a \textit{face component} of
$L$.  To describe the rest of $D$ which lies in $g$, we construct a
continuous map $\zv\co C\to \partial P$ (which is independent of $g$) such
that 1) $\zv$ maps vertices to vertices in the canonical way, 2) the
restriction of $\zv$ to every edge of $C$ is a homeomorphism onto the
canonically corresponding edge of $P$ and 3) the restriction of $\zv$ to
the inverse image of the interior of every face of $P$ is a homeomorphism.
The face $g$ of $C$ corresponds to two paired faces $f$ and $f^{-1}$ of
$P$.  Let $c$ be an edge cone of $f$ at an edge $e$ (as defined in the
fifth paragraph of Section~\ref{sec:generalizing}).  The face-pairing
$\ze$ pairs $c$ with an edge cone $c'$ of $f^{-1}$ at an edge $e'$.  Then
part of one component of $L$ has a projection $\zb$ in
$g\setminus\{\infty\}$ such that 1) only the endpoints of $\zb$ lie in an
edge of $g$, 2) $\zv(\zb)$ begins at the barycenter of $e$, 3) then an
initial segment of $\zv(\zb)$ lies in $c$, 4) then $\zb$ crosses under
$\za$, 5) then $\zb$ crosses over $\za$, 6) then a terminal segment of
$\zv(\zb)$ lies in $c'$ and 7) finally $\zv(\zb)$ ends at the barycenter
of $e'$.  The corridor complex $C$ is constructed so that we may, and do,
choose the projections $\zb$ for a fixed $g$ (and $f$) and varying $c$ to
have no self-crossings and no crossings with each other.
Figure~\ref{fig:corarcs} shows the face component $\za$ and the
projections $\zb$ for the face on the right part of
Figure~\ref{fig:corridor}. Constructing
such projections for every face $g$ of $C$ obtains $D$.  The components of
$L$ other than the face components are in bijective correspondence with
the edge cycles of $\ze$.  We call these components of $L$ \textit{edge
components}.  We call $D$ a \textit{corridor complex link diagram}
for $\ze$. We call $L$ a \textit{corridor complex link} for $\ze$.

\begin{figure}[ht!]
\centerline{\includegraphics{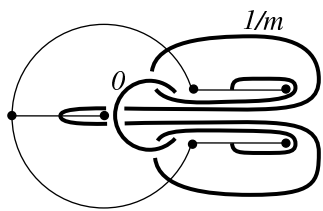}} \caption{A framed corridor
complex link diagram for Example~\ref{ex:s2xs1}} \label{fig:s2xs1d}
\end{figure}

\begin{ex}[]\label{ex:s2xs1'} We return to the model
face-pairing in Example~\ref{ex:s2xs1}.  A corridor complex for it
appears in Figure~\ref{fig:s2xs1d}, drawn with thin arcs.  A framed
corridor complex link diagram for it also appears in
Figure~\ref{fig:s2xs1d}, drawn with thick arcs.  The model face-pairing
has only one edge cycle, and we let it have multiplier $m$.
Theorem~\ref{thm:link} states that the associated twisted face-pairing
manifold $M$ is obtained by Dehn surgery on the framed link in
Figure~\ref{fig:s2xs1d}.  The framed link in
Figure~\ref{fig:s2xs1d} is isotopic to a link consisting of two unlinked
circles, one with framing 0 and one with framing $1/m$.  As in
Proposition 14.4 of \cite{PS}, Dehn surgery on a circle in $S^3$ with
framing 0 gives $S^2\times S^1$, and as in Proposition 14.6 of
\cite{PS}, Dehn surgery on a circle in $S^3$ with framing $1/m$ gives
$S^3$.  Thus $M$ is the connected sum of $S^2\times S^1$ and $S^3$.  In
other words, $M$ is $S^2\times S^1$ for every choice of the multiplier
$m$.
\end{ex}

\begin{thm}[]\label{thm:link}
Let $M=M(\ze,\text{mul})$ be a twisted
face-pairing manifold, and let $E_1,\dotsc,E_m$ be the edge cycles of
$\ze$. Let $L$ be a corridor complex for $\ze$ with diagram $D$. Define a
framing on $L$ as follows: every face component of
$L$ has framing 0, and the edge component of $L$ corresponding to
$E_i$ has framing $\text{mul}(E_i)^{-1}$ plus its blackboard framing
relative to $D$ for every $i\in\{1,\dotsc,m\}$.  Then the manifold
obtained by Dehn surgery on the framed link $L$ is homeomorphic to $M$.
\end{thm}
\begin{proof} Let $C$ be the corridor complex for $\ze$ from
which $D$ is constructed.  As in the construction of $D$, we view the
underlying space of $C$ as the one-point compactification $\bR^2\cup
\{\infty\}$ of $\bR^2$, where the point $\infty$ lies in the interior of
some face of $C$.  We choose standard coordinates $x$, $y$ and $z$ for
$\bR^3$, and we identify $C\setminus\{\infty\}$ with the $xy$-plane in
$\bR^3$.  We choose a closed standard metric ball in $\bR^3$ centered at
the origin so large that it contains every edge of $C$ in its interior.
Let $X$ be the solid hemisphere consisting of all points of this ball on
and below the $xy$-plane.

In this paragraph we construct a handlebody in $\bR^3$ by attaching
handles to $X$.  Let $f$ and $f^{-1}$ be two paired faces of $P$.  Let
$g$ be the face of $C$ corresponding to $f$ and $f^{-1}$.  If
$\infty\notin g$, then $g\subseteq \partial X$.
If $\infty\in g$, then $g\cap \partial X$
has  nonempty interior and is the complement in $g$ of a
neighborhood of $\infty$.  We attach a standard handle to $g\cap X$.
This handle is embedded in $\bR^3$ so that its vertical projection to
the $xy$-plane lies both in $X$ and in the interior of $g$.
Figure~\ref{fig:handle1} gives a view from above of $g$ and the handle
attached to $g$, where both $f$ and $f^{-1}$ are squares joined by a
simple corridor.  Figure~\ref{fig:handle2} gives another view of this
handle.  For every two paired faces of $P$ we attach a handle to $X$ in
this way.  We denote the result by $H$.  It is clear that $H$ is a
handlebody and that the closure of the complement of $H$ in $S^3$ is
also a handlebody.

\begin{figure}[ht!]
\centerline{\includegraphics{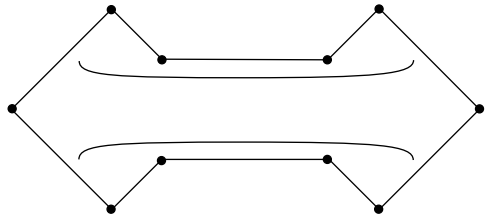}} \caption{Top view of the handle
attached to $g$} \label{fig:handle1}
\end{figure}

\begin{figure}[ht!]
\centerline{\includegraphics{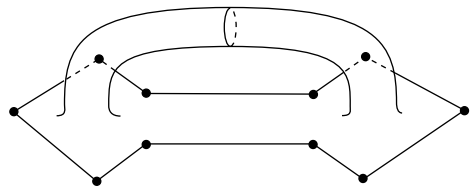}} \caption{Another view of the
handle attached to $g$} \label{fig:handle2}
\end{figure}

We next construct simple closed curves in $\partial H$ as follows.
First choose a barycenter for every edge of $C$.  Again let $f$ and
$f^{-1}$ be two paired faces of $P$, and let $g$ be the corresponding
face of $C$.  Just as in the construction of $D$, construct curves in
$\partial H$ which lie in and above $g$; these curves cross the handle
and they join barycenters of edges of $g$ which correspond to edges of
$f$ and barycenters of edges of $g$ which correspond to edges of
$f^{-1}$.  For every corridor edge $e$ of $g$ construct an arc in $g$ in
the obvious way which joins the barycenter of $e$ and the barycenter of
the edge of $g$ across the corridor from $e$.  We construct all these
curves so that only their endpoints lie in edges of $g$ and they are
pairwise disjoint except possibly at endpoints.  Finally, construct a
meridian curve for the handle of $H$ attached to $g$ such that this
meridian curve meets each of the curves which cross the handle exactly
once.  Figure~\ref{fig:handle3} shows a top view of $g$ and the handle
of $H$ attached to $g$ with the curves just constructed drawn with thick
solid and dashed arcs.  Doing this for every two paired faces of $P$, we
obtain two families of simple closed curves in $\partial H$.  The curves
$\zg_1,\dotsc,\zg_n$ in one family are the meridian curves of the
handles of $H$.  The curves $\zd_1,\dotsc,\zd_m$ in the other family
correspond canonically to the edge cycles of $\ze$.

\begin{figure}[ht!]
\centerline{\includegraphics{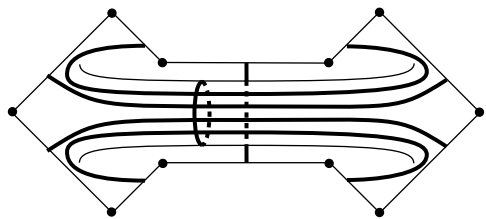}} \caption{Constructing curves in
the part of $\partial H$ which project vertically to $g$}
\label{fig:handle3}
\end{figure}

Let $S$ be the edge pairing surface of the twisted face-pairing $\zd$,
let $\za_1,\dotsc,\za_n$ be the vertical meridian curves, and let
$\zb_1,\dotsc,\zb_m$ be core curves of the edge cycle cylinders.
Then $S$ is homeomorphic to the edge pairing surface $S'$ of $\ze$ by a
homeomorphism that takes $\za_1,\dotsc,\za_n$ to the vertical meridian
curves of $S'$ and takes $\zb_1,\dotsc,\zb_m$ to the diagonal meridian
curves of $S'$. The surface $\partial H$ is homeomorphic to the edge
pairing surface $S'$ of $\ze$ by a
homeomorphism that takes $\zg_1,\dotsc,\zg_n$ to the vertical meridian
curves of $S'$ and takes $\zd_1,\dotsc,\zd_m$ to the diagonal meridian
curves of $S'$. Hence
the curves $\zg_1,\dotsc,\zg_n$ and $\zd_1,\dotsc,\zd_m$ can
be indexed so that there exists a homeomorphism $\zv\co S\to \partial H$
such that $\zv(\za_i)=\zg_i$ and $\zv(\zb_j)=\zd_j$ for every
$i\in\{1,\dotsc,n\}$ and $j\in\{1,\dotsc,m\}$.  Theorem~\ref{thm:surgery}
produces a framed link $L$ in $S^3$ such that the manifold obtained by
Dehn surgery on $L$ is homeomorphic to $M$.  Finally, it is clear that $D$
is a diagram of $L$ and that the framings are as claimed.

This proves Theorem~\ref{thm:link}.
\end{proof}

\section{Examples}\label{sec:examples}\nosubsections

In this section we present some examples in which we use
Theorem~\ref{thm:link} to identify some twisted face-pairing manifolds.
We have already given such an example in Example~\ref{ex:s2xs1'}, where
we constructed a framed link in $S^3$ for the model face-pairing in
Example~\ref{ex:s2xs1}.  Using this we showed that the twisted
face-pairing manifolds in Example~\ref{ex:s2xs1} are all homeomorphic to
$S^2\times S^1$.

\begin{figure}[ht!]
\centerline{\includegraphics{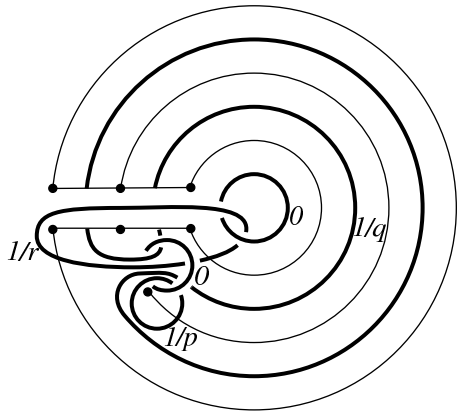}} \caption{The corridor complex
and framed link diagram for Example~\ref{ex:target}} \label{fig:targetd}
\end{figure}

\begin{figure}[ht!]
\centerline{\includegraphics{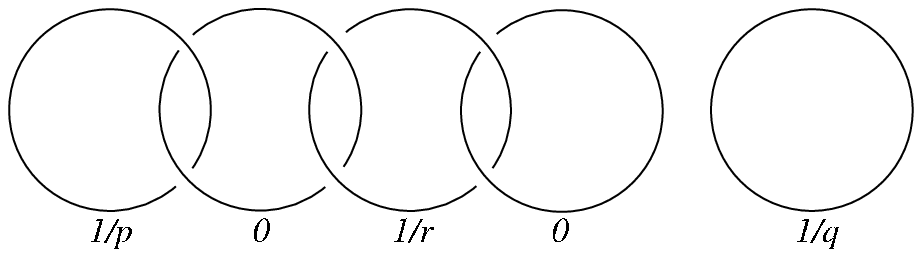}} \caption{A simpler framed
link} \label{fig:targetd1}
\end{figure}

\begin{figure}[ht!]
\centerline{\includegraphics{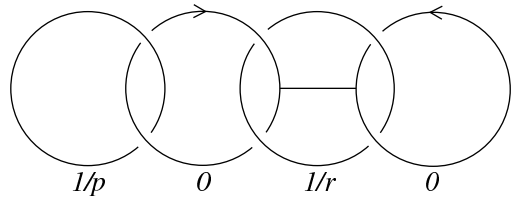}} \caption{Preparing for a Kirby
move of type 2} \label{fig:targetd2}
\end{figure}

\begin{figure}[ht!]
\centerline{\includegraphics{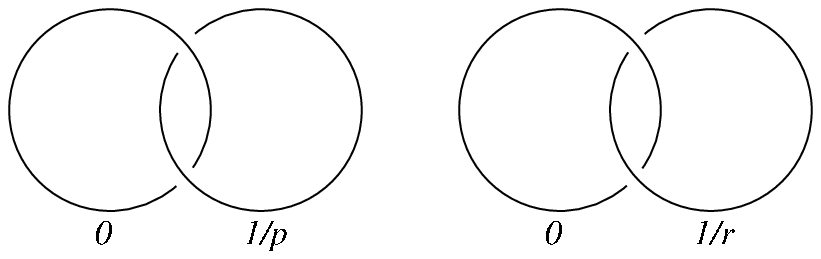}} \caption{A simpler framed
link} \label{fig:targetd3}
\end{figure}

\begin{ex}[]\label{ex:target''} We return to the model
face-pairing in Example~\ref{ex:target}.  We choose multipliers of the
edge cycles in line~\ref{lin:targetedge} to be $p$, $q$, and $r$, in
order.  A corridor complex for Example~\ref{ex:target} appears in
Figure~\ref{fig:targetd}, drawn with thin arcs.  A framed link diagram for
it also appears in Figure~\ref{fig:targetd}, drawn with thick arcs.
Theorem~\ref{thm:link} states that the associated twisted face-pairing
manifold $M$ is obtained by Dehn surgery on the framed link in
Figure~\ref{fig:targetd}.  The framed link in
Figure~\ref{fig:targetd} is isotopic to the framed link in
Figure~\ref{fig:targetd1}.  The component of the link in
Figure~\ref{fig:targetd1} with framing $1/q$ corresponds to a connected
summand of $M$.  But, as in Example~\ref{ex:s2xs1'}, this connected
summand is $S^3$.  So we delete the component of the link in
Figure~\ref{fig:targetd1} with framing $1/q$.  We modify the component of
the link in Figure~\ref{fig:targetd1} with framing 0 which links the
components with framings $1/p$ and $1/r$ by means of a Kirby move of type
2.  For this we orient the components with framing 0 and connect them with
an arc as shown in Figure~\ref{fig:targetd2}.  The result is a link
isotopic to the one in Figure~\ref{fig:targetd3}.  It follows from
Proposition 17.3 of \cite{PS} that $M$ is a connected sum of the lens
space $L(p,-1)=L(p,1)$ and the lens space $L(r,-1)=L(r,1)$.
\end{ex}

\begin{figure}[ht!]
\centerline{\includegraphics{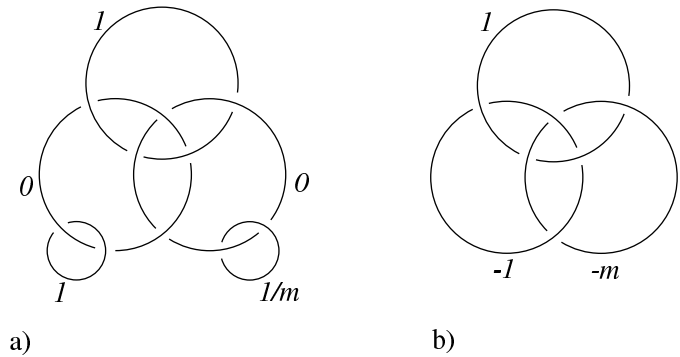}} \caption{Two framed
links for Example~\ref{ex:tet1a}} \label{fig:tet1ad12}
\end{figure}

\begin{figure}[ht!]
\centerline{\includegraphics{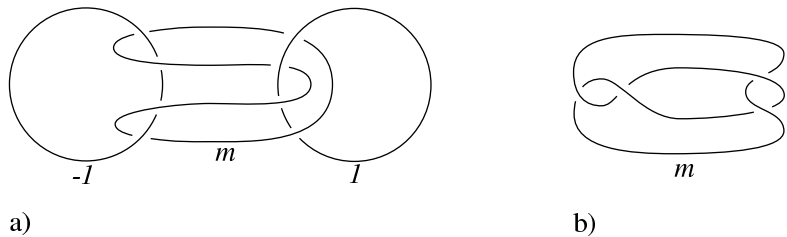}} \caption{Two more framed
links for Example~\ref{ex:tet1a}}
\label{fig:tet1ad34}
\end{figure}

\begin{ex}[]\label{ex:tet1a} We return to the model face-pairing which we
described at the beginning of the introduction.
We choose multipliers $m_1=1$, $m_2=1$, and $m_3=m$. A
corridor complex for this example appears in Figure~\ref{fig:tet1pdi3},
drawn with thin arcs, and a framed link diagram for it also appears in
Figure~\ref{fig:tet1pdi3}, drawn with thick arcs. The
part of the link in Figure~\ref{fig:tet1pdi3} which is the union of the
components with framing 0 and the component which in the diagram crosses
both components with framing 0 is isotopic to the Borromean rings.  So the
framed link in Figure~\ref{fig:tet1pdi3} is isotopic to the link in
part a) of Figure~\ref{fig:tet1ad12}.  We simplify the framed link in
part a) of Figure~\ref{fig:tet1ad12} using Kirby calculus by performing twist
moves, which are discussed in Sections 16.4, 16.5 and 19.4 of \cite{PS}
under the name Fenn-Rourke moves.
Twisting $-m$ times along the component with framing $1/m$, twisting $-1$
times along the similar component with framing 1, and deleting resulting
components with framing $\infty$ yields the link in
part b) of Figure~\ref{fig:tet1ad12}.
Because the link in part b) of Figure~\ref{fig:tet1ad12} is
amphicheiral we may, and do, multiply every framing by $-1$.  We isotop
the result to the framed link in part a) of Figure~\ref{fig:tet1ad34}.
Now we perform
twist moves on the link in part a) of Figure~\ref{fig:tet1ad34}.
We twist 1
time along the component with framing $-1$, twist $-1$ times along the
component with framing 1, and delete resulting components with framing
$\infty$.  The result is shown in part b) of Figure~\ref{fig:tet1ad34}.
This is the
figure eight knot with framing $m$.  If $m=1$, then $M$ is the Brieskorn
homology sphere $\zS(2,3,7)$, which has the geometry of the universal
cover of $PSL(2,\bR)$. According to Theorem 4.7 of \cite{T}, $M$ is
hyperbolic if $m\ge 5$.
\end{ex}

\begin{figure}[ht!]
\centerline{\includegraphics{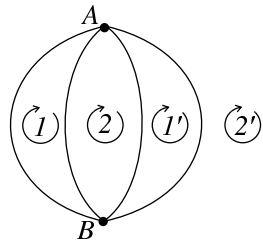}} \caption{The complex $P$ for
Example~\ref{ex:nil}} \label{fig:nilp}
\end{figure}

\begin{figure}[ht!]
\centerline{\includegraphics{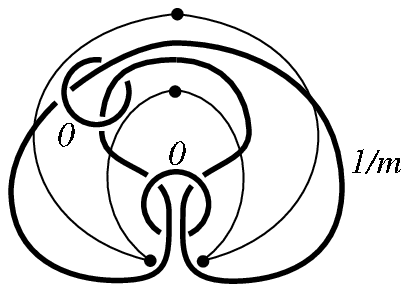}} \caption{A corridor complex and
framed link diagram for Example~\ref{ex:nil}} \label{fig:nild}
\end{figure}

\begin{figure}[ht!]
\centerline{\includegraphics{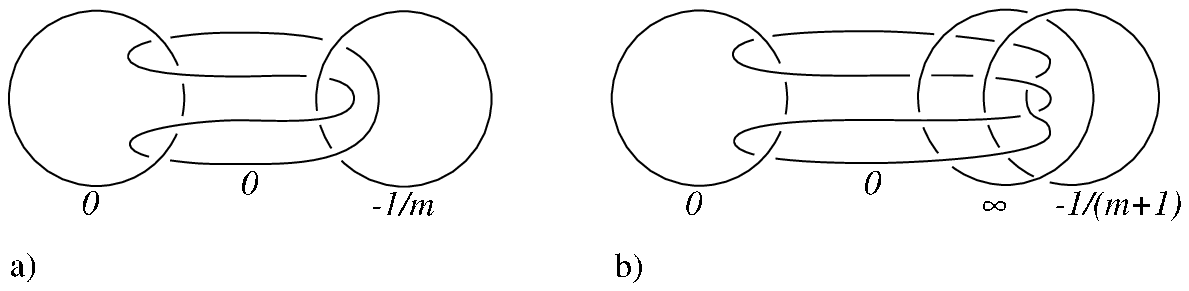}} \caption{Two framed links
for Example~\ref{ex:nil}}
\label{fig:nild12}
\end{figure}

\begin{figure}[ht!]
\centerline{\includegraphics{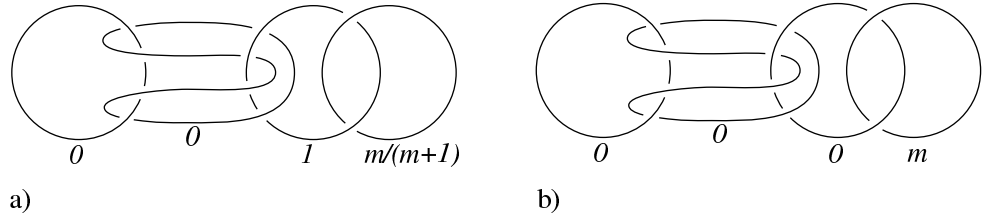}} \caption{Two more framed
links for Example~\ref{ex:nil}}
\label{fig:nild34}
\end{figure}

\begin{ex}[]\label{ex:nil} This example is closely related to the
previous one.  The model faceted 3-ball for this example is gotten from
the faceted 3-ball given in Figure~\ref{fig:tet1p} by collapsing the
edge $AB$ to a point and collapsing the edge $CD$ to a point.  The
result is the faceted 3-ball $P$ given in Figure~\ref{fig:nilp}.
Because the edges $AB$ and $CD$ in Figure~\ref{fig:tet1p} are both
fixed by the model face-pairing of Example~\ref{ex:tet1a}, the model
face-pairing of Example~\ref{ex:tet1a} induces a model face-pairing
$\ze$ on $P$.  The face-pairing $\ze$ pairs the faces of $P$ as
indicated in Figure~\ref{fig:nilp}, and the face-pairing maps of $\ze$
fix the vertices $A$ and $B$.  The model face-pairing $\ze$ has one edge
cycle.  This edge cycle has length 4 and corresponds to the edge
cycle of length 4 in Example~\ref{ex:tet1a}.
We let this edge cycle of $\ze$ have
multiplier $m$.  A corridor complex for $\ze$ appears in
Figure~\ref{fig:nild}, drawn with thin arcs, and a framed link diagram
for it also appears in Figure~\ref{fig:nild}, drawn with thick arcs.
This link is the Borromean rings.  As in Example~\ref{ex:tet1a} we may,
and do, multiply the framings by $-1$ and we isotop the link in
Figure~\ref{fig:nild} to obtain the framed link in
part a) of Figure~\ref{fig:nild12}.  Now we perform a twist move by twisting
$-1$ times along the component with framing $-1/m$ and we introduce a
component with framing $\infty$ to obtain the framed link in part b) of
Figure~\ref{fig:nild12}.  Next we twist 1 time along the component with
framing $\infty$ to obtain the link in
part a) of Figure~\ref{fig:nild34}.  Finally,
we twist $-1$ times along the component with framing $m/(m+1)$ to obtain
the link in part b) of Figure~\ref{fig:nild34}.
Adding a component with framing $\infty$ that is parallel to the
component with framing $m$ gives a
special case of the link at the top of Figure 12 on page 146 of
\cite{M}.  It follows that $M$ is the Seifert fibered manifold $\langle
Oo1|0;(m,1)\rangle$ in the notation of \cite{M}.  This means that $M$ is
orientable with an orientable base surface of genus 1, that the Euler
number of $M$ is 0, and that $M$ has one exceptional fiber of type
$(m,1)$.  When $m=1$, the manifold $M$ is the Heisenberg manifold, the
prototype for Nil geometry.
\end{ex}

\np

\Addresses\recd
\end{document}

%% file: agtout.tex
\def\ifplaintex{\expandafter\ifx\csname documentclass\endcsname\relax}

\def\gtp{{\mathsurround=0pt\it $\cal G\mskip-2mu$eometry \&\ 
$\cal T\!\!$opology $\cal P\!$ublications}}  

\def\recd{{\small Received:\qua\receiveddate\ifx\reviseddate\relax
\else\qquad Revised:\qua\reviseddate\fi\par}} 

\def\lognumber#1{\def\thelognumber{#1}}
\def\volumenumber#1{\def\thevolumenumber{#1}}
\def\volumeyear#1{\def\thevolumeyear{#1}}
\def\papernumber#1{\def\thepapernumber{#1}}
\def\pagenumbers#1#2{\def\startpage{#1}\def\finishpage{#2}}
\def\published#1{\def\publishdate{#1}}

\def\received#1{\def\receiveddate{#1}}
\def\revised#1{\def\reviseddate{#1}}
\def\accepted#1{\def\accepteddate{#1}}

\def\asciiaddress#1{\def\theasciiaddress{#1}}

\let\\\par\let\thelognumber\relax\let\thevolumenumber\relax
\let\thepapernumber\relax\let\thevolumeyear\relax\let\startpage\relax
\let\finishpage\relax\let\publishdate\relax\let\receiveddate\relax
\let\reviseddate\relax\let\accepteddate\relax\let\theasciititle\relax
\let\theasciiauthors\relax\let\theasciiaddress\relax
\let\theasciiabstract\relax

\let\theasciiemail\relax

\ifplaintex
\font\logobig=cmssbx10 scaled 3836
\font\logomed=cmssbx10 scaled 2557
\else
\font\logobig=cmssbx10 scaled 4200
\font\logomed=cmssbx10 scaled 2800
\fi

\long\def\makeagttitle{   
\count0=\startpage
\agt\hfill      
\hbox to 45truept{\vbox to 0pt{\vglue -13truept{\logomed A\kern -.37em{\logobig 
T}\kern -.38em G}\vss}\hss}
\break
{\small Volume \thevolumenumber\ (\thevolumeyear)
\startpage--\finishpage\nl
Published: \publishdate}

\vglue .25truein

{\parskip=0pt\leftskip 0pt plus
1fil\def\\{\par\smallskip}{\Large\bf\thetitle}\par\medskip} \vglue
0.05truein

{\parskip=0pt\leftskip 0pt plus 1fil\def\\{\par}{\sc\theauthors}
\par\medskip}%
 
\vglue 0.03truein 

{\small\leftskip 25truept\rightskip 25truept{\bf Abstract}\stdspace\theabstract

{\bf AMS Classification}\stdspace\theprimaryclass
\ifx\thesecondaryclass\relax\else; \thesecondaryclass\fi\par
{\bf Keywords}\stdspace \thekeywords\par}\vglue 7truept

}   

\ifplaintex
\font\phead=cmsl9 scaled 950
\font\pnum=cmbx10 scaled 913
\font\pfoot=cmsl9 scaled 950
\headline{\vbox to 0pt{\vskip -4.5mm\line{\small\phead\ifnum
\count0=\startpage ISSN 1472-2739 (on-line) 1472-2747 (printed)
\hfill {\pnum\folio}\else\ifodd\count0\def\\{ }%
\ifx\theshorttitle\relax\thetitle\else\theshorttitle\fi\hfill{\pnum\folio}
\else\def\\{ and }{\pnum\folio}\hfill\ifx\theshortauthors\relax\theauthors
\else\theshortauthors\fi\fi\fi}\vss}}
\footline{\vbox to 0pt{\vglue 0mm\line{\small\pfoot\ifnum\count0=\startpage
\copyright\ \gtp\hfill\else
\agt, Volume \thevolumenumber\ (\thevolumeyear)\hfill\fi}\vss}}
\else
\font=cmsl9 scaled 1050
\font\lnum=cmbx10 
\font=cmsl9 scaled 1050
\makeatletter
\def\@oddhead{{\small\lhead\ifnum\count0=\startpage ISSN 1472-2739 
(on-line) 1472-2747 (printed)\hfill {\lnum\number\count0}\else\ifodd\count0
\def\\{ }\ifx\theshorttitle\relax \thetitle \else\theshorttitle\fi\hfill
{\lnum\number\count0}\else\def\\{ and }{\lnum\number\count0}
\hfill\ifx\theshortauthors\relax 
\theauthors\else\theshortauthors\fi\fi\fi}}\def\@evenhead{\@oddhead}
\def\@oddfoot{\small\lfoot\ifnum\count0=\startpage\copyright\ \gtp\hfill\else
\agt, Volume \thevolumenumber\ (\thevolumeyear)\hfill\fi}
\def\@evenfoot{\@oddfoot}
\makeatother
\fi
\let\maketitlepage\makeagttitle
\let\makeshorttitle\maketitlepage
\let\maketitle\maketitlepage

\newwrite\gtoutfile
\long\gdef\makeheadfile{  
{\def\\{, }\def\s{ }
\immediate\openout\gtoutfile head.xxx
\immediate\write\gtoutfile{To: math@arxiv.org}
\immediate\write\gtoutfile{Subject: put OR rep NNNNN:ppppp}
\immediate\write\gtoutfile{--text follows this line--}
\immediate\write\gtoutfile{Proxy-for: \ifx\theasciiauthors\relax
\theauthors\else\theasciiauthors\fi\s<\ifx\theasciiemail\relax\theemail\else\theasciiemail\fi>}
\immediate\write\gtoutfile{\noexpand\\}
\immediate\write\gtoutfile{Authors: \ifx\theasciiauthors\relax
\theauthors\else\theasciiauthors\fi}
{\def\\{ }\immediate\write\gtoutfile{Title: \ifx\theasciititle\relax
\thetitle\else\theasciititle\fi}}
\immediate\write\gtoutfile{Subj-class: GT or SG, GR etc}
\immediate\write\gtoutfile{MSC-class: \theprimaryclass\ifx\thesecondaryclass\relax\else, \thesecondaryclass\fi}
\immediate\write\gtoutfile{Journal-ref: Algebr. Geom. Topol. \thevolumenumber\s
(\thevolumeyear) \startpage-\finishpage}
\immediate\write\gtoutfile{Comments: Published by Algebraic and
Geometric Topology at}
\immediate\write\gtoutfile{\s\s\s  http://www.maths.warwick.ac.uk/agt/AGTVol\thevolumenumber/agt-\thevolumenumber-\thepapernumber.abs.html}
\immediate\write\gtoutfile{\noexpand\\}
\immediate\write\gtoutfile{}
\ifx\theasciiabstract\relax
\immediate\write\gtoutfile{\theabstract}\else
\immediate\write\gtoutfile{\theasciiabstract}\fi
\immediate\write\gtoutfile{}
\immediate\write\gtoutfile{\noexpand\\}
\immediate\write\gtoutfile{}
\immediate\closeout\gtoutfile}}  

\def\maketitlepage{\makeagttitle\makeheadfile}
\let\makeshorttitle\maketitlepage
\let\maketitle\maketitlepage